\documentclass{article}
\usepackage{arxiv}

\usepackage{authblk}
\usepackage{enumerate} 
\usepackage[utf8]{inputenc}
\usepackage[normalem]{ulem}
\usepackage{url}
\usepackage{eurosym}
\usepackage{bm}
\usepackage{subfigure}
\usepackage{amsmath}
\usepackage{multirow}
\usepackage{booktabs}
\usepackage{caption}
\usepackage{svg}
\usepackage{graphicx}
\usepackage{hyperref}
\usepackage{float}
\usepackage[version=4]{mhchem}
\usepackage{bbold}
\usepackage{longtable,tabularx}
\usepackage{multirow}
\usepackage{amsmath}
\usepackage{algorithm, setspace}
\usepackage{algpseudocode}
\usepackage{siunitx,letltxmacro}
\LetLtxMacro{\svqty}{\qty}
\usepackage{physics}
\usepackage{comment}
\LetLtxMacro{\qty}{\svqty}
\usepackage[colorinlistoftodos]{todonotes}
\usepackage[nopostdot,acronym,nogroupskip,nonumberlist]{glossaries}

\definecolor{myred}{RGB}{192, 0, 0}
\definecolor{mygreen}{RGB}{152,185,84}
\definecolor{myblue}{RGB}{0, 70, 127}

\usepackage{tikz}

\usepackage{float}	 											

\usepackage{setspace} 											
\usepackage{color}
\usepackage{caption}											
\newcommand\newl[2]{ 											
 	\expandafter\newlength\csname #1\endcsname
 	\expandafter\setlength\csname #1\endcsname{#2}}
\newcommand\setl[2]{ 											
 	\expandafter\setlength\csname #1\endcsname{#2}}

\newl{figh}{.0\columnwidth}										
\newl{figw}{.0\columnwidth}										
\newl{boxh}{.0\columnwidth}										
\newl{boxw}{.0\columnwidth}										
\newl{tcw}{.0\columnwidth}										

\newl{subfigskip}{.6em}											
\newl{plotvsep}{.8em}											

\usepackage{pgfplots}
\usepgfplotslibrary{groupplots}
\usetikzlibrary{positioning}

\pgfplotsset{
	compat=newest,
	plot coordinates/math parser=false,
}

\pgfdeclarelayer{annotation}
\pgfdeclarelayer{grid}
	\pgfsetlayers{grid,main,annotation}
\newcommand{\addimagegrid}{1}									
\newcommand{\imagegrid}{										
	\begin{pgfonlayer}{grid}
		\draw[black!5, ultra thin, step=1,] (img.south west) grid (img.north east);
	 	\draw[black!15, ultra thin, step=10,] (img.south west) grid (img.north east);
	 	\foreach \x in {0,10,...,100} { \node [below] at (\x,0) {\x};}
	 	\foreach \y in {0,10,...,100} { \node [left] at (0,\y) {\y};}
	\end{pgfonlayer}
}		

\usepackage{amssymb}

\begin{document}

\title{Team theAntipodes: Solution Methodology for GTOC11}

\author[$\dagger$]{Roberto Armellin}
\author[$\star$]{Laurent Beauregard}
\author[$\ddag$]{Andrea Bellome}
\author[$+$]{Nicolò Bernardini}
\author[$\diamond$]{Alberto Foss\`a}
\author[$+$]{Xiaoyu Fu}
\author[$+$]{Harry Holt}
\author[$\circledast$]{Cristina Parigini}
\author[$\dagger$]{Laura Pirovano}
\author[$\dagger$]{Minduli Wijayatunga}

\affil[$\dagger$]{Te P\=unaha \=Atea - Space Institute, University of Auckland, 1010, Auckland, New Zealand}
\affil[$\star$]{Mission Analysis Section, Telespazio@ESA/ESOC, Robert-Bosch Strasse 5, 64293 Darmstadt, Germany}
\affil[$\ddag$]{Space Research Group, Cranfield University, MK43 0AL, Cranfield, United Kingdom}
\affil[$+$]{Surrey Space Centre, University of Surrey, GU2 7XH, Guildford, United Kingdom}
\affil[$\diamond$]{Institut Sup\'erieur de l’A\'eronautique et de l’Espace, 31055, Toulouse, France}
\affil[$\circledast$]{School of Mathematical Sciences, University of Southampton, SO17 1BJ, Southampton, United Kingdom}

\maketitle

\begin{abstract}

This paper presents the solution approach developed by the team \lq\lq theAntipodes\rq\rq\ for the 11th Global Trajectory Optimization Competition (GTOC11).
The approach consists of four main blocks: 1) mothership chain generation, 2) rendezvous table generation, 3) the dispatcher, and 4) the refinement. Blocks 1 and 3 are purely combinatorial optimization problems that select the asteroids to visit and allocate them to the Dyson ring stations. The rendezvous table generation involves interpolating time-optimal transfers to find all transfer opportunities between selected asteroids and the ring stations. The dispatcher uses the data stored in the table and allocates the asteroids to the Dyson ring stations in an optimal fashion. The refinement ensures each rendezvous trajectory meets the problem accuracy constraints, and introduces deep-space maneuvers to the mothership transfers. We provide the details of our solution that, with a score of 5,992, was worth 3rd place.      

\end{abstract}

\section{Introduction} \label{introduction}

The 11th Global Trajectory Optimization Competition (GTOC11) considered the problem of designing a \lq\lq Dyson Ring" using asteroids in the solar system. A summary of the problem is provided below, the interested reader can find all the details on the competition web-page \url{gtoc11.nudt.edu.cn}. 

\begin{figure}[!ht]
\begin{center}
\includegraphics[width = \textwidth]{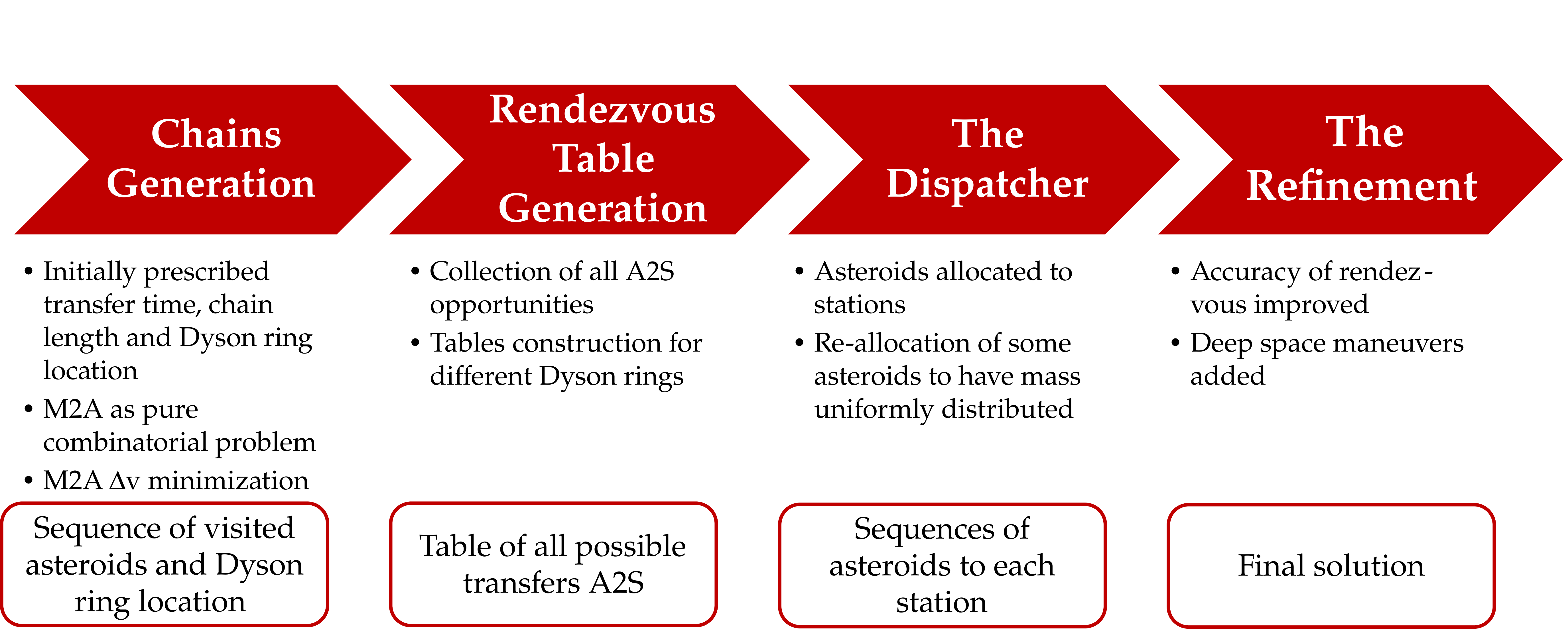}
\captionsetup{justification=centering}
\caption{Flowchart representing the solution's computational blocks. Boxes at the bottom represent the block's output, which served as input for the next block. Transfers performed by the motherships to visit asteroids are shortened to M2A, while transfers from asteroids to the stations as A2S.}
\label{fig:flowchart}
\end{center}
\end{figure}

Up to 10 motherships can be launched from the Earth to perform flybys with  asteroids selected from a list of approximately 83,000 objects. These transfers are realized with impulsive maneuvers. A maximum of four impulses can be used by a mothership in each transfer leg, and an impulse of 6 km/s is available at Earth-departure for free. At an asteroid encounter (i.e., when the mothership and an asteroid are within 1 km with a relative velocity less than 2 km/s) a multi-functional asteroid transfer device (ATD) is place on the asteroid. Once activated, the ATD consumes the asteroid resources to produce thrust, enabling its transfer to a designated circular Dyson ring orbit. These asteroids are used to build 12 stations sequentially to generate solar power. These stations are distributed uniformly in phase. The semi-major axis ($a_{D}$), inclination ($i_{D}$), right ascension of ascending node ($\Omega_{D}$), and phase of the first station ($\phi_{S_1}$) are design parameters. The objective is to maximize the function
\begin{equation}\label{eq:obj}
J_{GTOC} = B \frac{10^{-10} \cdot m_{\mathrm{min}}^{\text S}}{a_{\text {D }}^{2} \sum_{i=1}^{10}\left(1+\Delta v_{i} / 50\right)^{2}},
\end{equation}
in which $B$ is a submission date-based reward, $m_{\mathrm{min}}^{\text {S }}$ is the mass of the least massive station, and $\Delta v_{i}$ is the total $\Delta v$ of the $i$-th mothership. Other relevant problem parameters and constraints are:
\begin{enumerate}[i)]
    \item All the events must occur between January 1, 2121 00:00:00 UT and January 1, 2141 00:00:00 UT, inclusive.
    \item Once placed on an asteroid, an ATD takes a minimum of 30 days before it can be activated.
    \item An asteroid obtains a fixed acceleration of 10$^{-4}$ m/s$^2$ once the ATD is activated. This acceleration can only be switched off at ring arrival. 
    \item  An asteroid's mass decreases as $\dot{m}^{A} = -\alpha m^{A}_0$, with $\alpha = 6\cdot 10^{-9}$ s$^{-1}$ and $m^{A}_0$ the initial mass of the asteroid.
    \item \label{constraint90days} The time interval between the last asteroid arriving at one station and the first asteroid arriving at the next station cannot be shorter than 90 days.
    \item The distance of motherships and asteroids from the Sun cannot get smaller than 0.4 AU. 
    \item The minimum value of $a_{D}$ is 0.65 AU.
\end{enumerate}
This paper summarizes the solution methodology developed by the team \lq\lq theAntipodes\rq\rq\ who submitted the third best solution with a score of $J = 5992$.
The approach is based on four main computational blocks: 1) mothership chains generation, 2) rendezvous table generation, 3) the dispatcher, and 4) the refinement, as shown in Fig.~\ref{fig:flowchart}. 
{With the exception of basic astrodynamics libraries (e.g., ephemeris, orbital elements conversion, Lambert's solver, Edelbaum's approximation) all the key computational blocks were developed from scratch by the team during the competition. Third-party nonlinear solvers and optimization routines were then used to compute numerical solutions to the formulated problems.}

The mothership chain generation is split into three parts. Firstly, a preliminary assessment of transfer times, chain length, and Dyson ring location is conducted using a genetic algorithm (GA) optimizer \cite{mitchell1998introduction} as described in Sec.~\ref{Sec:CAS}. To do so, Lambert's arcs \cite{izzo2015revisiting} are used to design the mothership transfers and {Edelbaum's analytic expressions (\cite{Kluever2011} Eq. 1--4)} are used to estimate the time for an asteroid to reach the ring. Next, this preliminary assessment informs the transcription of the mothership transfer design into a pure combinatorial problem, solved with a beam search (BS) algorithm in Sec.~\ref{Sec:CO}. The output of this process is the sequence of visited asteroids and a rough estimate of the encounter dates. Finally, for a fixed asteroid sequence, each {mothership} $\Delta v$ is minimized by optimizing {all} the encounter dates and pre- and post-flyby impulses (see Sec.~\ref{sec:R}). At this stage the Dyson ring parameters are also updated. 

Once the motherships transfers are computed, the focus shifts to the Dyson ring construction. Sec.~\ref{sec:TG} provides details on the rendezvous table generation. For all asteroids visited by the motherships, the table collects all minimum-time transfer opportunities to any station within the mission duration. To achieve this result, we first formulate and solve time-optimal orbit transfers using an indirect method \cite{rasotto2016multi}. For each asteroid, 16 optimal transfers were computed taking initial conditions uniformly spaced over one orbital period. These solutions and the corresponding arrival phases at the ring are interpolated to generate very accurate guesses for the rendezvous problem, thus populating the table, as described in Sec.~\ref{sec:ETF}. 

The rendezvous table provides the dispatcher input. The dispatcher determines the station building sequence and allocates asteroids to each station to maximize the objective function in Eq.~\eqref{eq:obj}. This combinatorial problem is formulated by dividing the objective function evaluation into two steps. Firstly, all the asteroids are preliminary allocated to the stations. Secondly, some asteroids are reallocated to distribute the mass uniformly among the stations, and, when convenient, some are dropped (as described in Sec.~\ref{sec:4.1}). The problem is initially solved with a GA followed by a more extensive search with a particle swarm optimizer (PSO) \cite{Eberhart1995}. For the same set of 10 motherships, the dispatcher is run on a set of rendezvous tables prepared for different $a_{D}$ in an attempt to find its optimal location.  

The refinement, detailed in Sec.~\ref{sec:5}, is the final step of our solution process. Firstly, the A2S transfers selected by the dispatcher are refined to achieve high-accuracy rendezvous. Secondly, mothership $\Delta v$'s are further reduced by enabling a deep-space maneuver in each transfer leg while fixing the flyby encounters. 

Lastly, Sec.~\ref{sec:OS} provides details of our optimal solution including trajectory visualizations. The solution output file can be downloaded from the competition web-page.    

\section{Mothership Chain Generation}\label{sec:2}

Designing missions that visit multiple orbital way-points is a notoriously challenging problem. There is an obvious need to manage the increasing complexity of the problem formulation, often through cleverly designed strategies that prune out the search space or manage the search strategy. They are generally transcribed in some form of global optimization, referred to as mixed-integer nonlinear programming (MINLP) problems \cite{schlueter2012nonlinear}. As such, the order of way-points to be visited is not known a priori but is instead the optimization objective. Thus, instances of MINLP usually involve the solution of a combinatorial optimization problem, coupled with optimal control theory.
Most approaches to solving the MINLP problem make use of three steps \cite{Alemany2007}: the first consists of defining a subset of potential targets (asteroids) to visit, based on their orbits and scientific characteristics. Secondly, a sequence of objects is found by means of global optimization algorithms. Lastly, the trajectories between asteroids are optimized with either local or global optimization. This scheme is also adopted by the authors here, as described in the following sections. 

\subsection{Beam Search Method} \label{first_step_transcription} 

With approximately 83,000 objects, the asteroid set presented for GTOC11 is already larger than any previous asteroid tour related GTOC competition problem. Exploring all possible $ N- $asteroids tours, for any $ N > 5$, represents an infeasible number to compute within a reasonable time window, as the number of asteroids sequences grows factorially with $ N$. Hence, the key to solving this problem must lie in managing this complexity efficiently. Thus, a reduction is performed on the data-set of available asteroids to enable a global search of transfers. {By inspecting the asteroids' orbital parameters and mass distribution, it was decided to exclude asteroids with $a_A > 2.8$ AU, $e_A>0.1584$, $i_A > 8.897$ deg, and $m_A< 5.8497 \times 10^{13}$ kg. This pruning resulted in a list of approximately 10,000 asteroids.}



Three fundamental assumptions are then made for the preliminary estimation of the key problem parameters. Firstly, the time of fight to reach the first asteroid, $\Delta t_{E2A}$, and the time of flight of each asteroid-to-asteroid (A2A) leg, $\Delta t_{A2A}$, are fixed and applied to all motherships. Secondly, the target Dyson ring is assumed to be on a circular orbit of radius $a_{D}$ with zero inclination. Thirdly, the orbital eccentricity of candidate asteroids is assumed to be zero for a preliminary analysis of low-thrust transfer to the ring. These choices were made predominantly to limit the aforementioned computational burden. Based on these assumptions, $\Delta t_{E2A}$, $\Delta t_{A2A}$, and $a_{D}$ are the three design parameters that enable the transcription of the problem into a pure combinatorial form.

\begin{figure}[!ht]
\begin{center}
\includegraphics[width=.5\columnwidth]{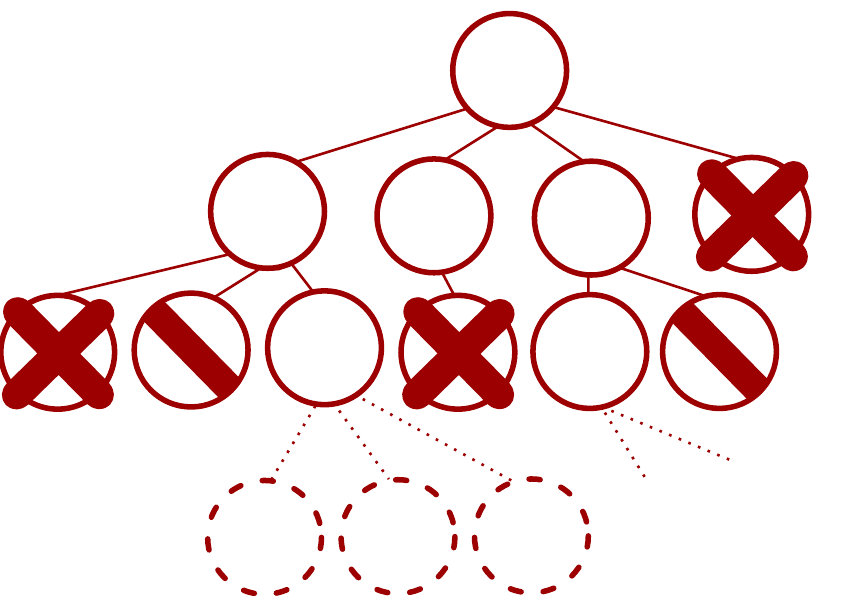}
\caption{Schematic of BS approach with incremental pruning. Barred nodes do not fall within the beam width and thus are pruned, while crossed nodes are pruned due to non-linear constraints violation. Dotted nodes are not explored yet.}
\label{fig:beam_search}
\end{center}
\end{figure}

A beam search (BS) strategy \cite{Shapiro1992, Wilt2010} with incremental pruning is employed to solve the combinatorial problem, as it has emerged as a standard approach in many past GTOC competitions. In BS algorithms, the computational effort is bounded by employing heuristics that prevent the exploration of non-promising branches. The BS algorithm is represented as a tree-graph in Fig.~\ref{fig:beam_search}. Each node encodes a trajectory with increasing number of asteroids per tree level. The exploration of possible trajectories is performed one depth-level at a time. From all the branches generated in one level, only a limited set, the so-called beam width (BW), are selected for expansion in successive nodes. {Specifically, at each level, the set of solutions are sorted with respect to mission-driven criteria (see Sec.~\ref{Sec:CAS}) and only the best (i.e., the BW) are kept for further expansion.} Selecting the proper BW is thus a compromise between solution quality and number, as well as computational effort. 

Branching, i.e. adding one trajectory leg, implies the solution of a Lambert problem between two consecutive asteroids for all the $ j-$nodes at a given tree level and $k-$available asteroids in the data-set. To prevent the tree expansion becoming intractable, the following pruning criteria are applied:

\begin{itemize}
    \item[1)] {Nodes exceeding with an Earth departure impulse exceeding 6 km/s are not further expanded.}
    \item[2)] Nodes with A2A cost $  \Delta v_{A2A} > 1.5 \mbox{ km/s} $ are not further expanded.
    \item[3)] Nodes with accumulated total $ \Delta v > 30 \mbox{ km/s} $ are not further expanded.
\end{itemize}

\begin{figure}[!ht]
\centering
	\begin{tikzpicture}
	\node[inner sep=0pt, anchor=south west] (img) at (0,0) {\includegraphics[width=0.4\textwidth]{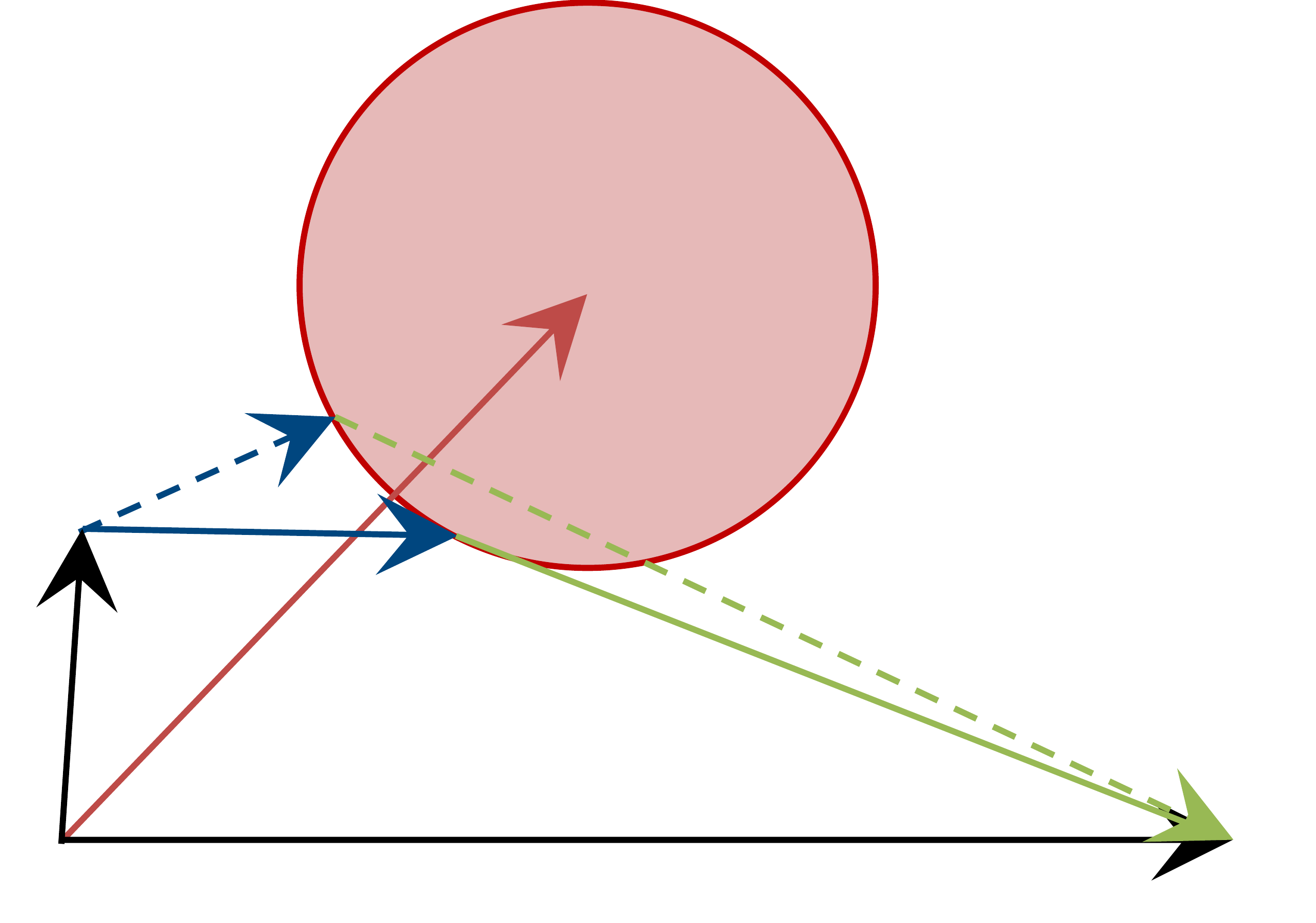}};
	\begin{scope}[ x={($0.01*(img.south east)$)}, y={($0.01*(img.north west)$)}, ]
		\ifnum\addimagegrid=0\imagegrid\fi
		\begin{pgfonlayer}{annotation}
 			\node[inner sep=0pt,anchor=north east] at (3,37) {$v^{-}$};
 			\node[inner sep=0pt,anchor=north east] at (89,4) {$v^{+}$};
 			\node[inner sep=0pt,anchor=north west] at (45,58) {{$v_{A}$}};
 			\node[inner sep=0pt,anchor=south east] at (15,50) {{$\Delta v_{1}$}};
 			\node[inner sep=0pt,anchor=south west] at (10,30) {{$\Delta v_{1}^{*}$}};
 			\node[inner sep=0pt,anchor=south west] at (64,31) {{$\Delta v_{2}$}};
 			\node[inner sep=0pt,anchor=north east] at (58,27) {{$\Delta v_{2}^{*}$}};
		\end{pgfonlayer}
	\end{scope}	
		\end{tikzpicture}
\caption{Velocity triangles at asteroid rendezvous. The non-starred impulses are the non-optimal one used in the mothership chain generation. The starred impulses are the optimal ones computed in Sec.~\ref{sec:R}. The ball around the asteroid velocity represents the 2 km/s constraints.}
\label{fig:velocity_triangles}
\end{figure}
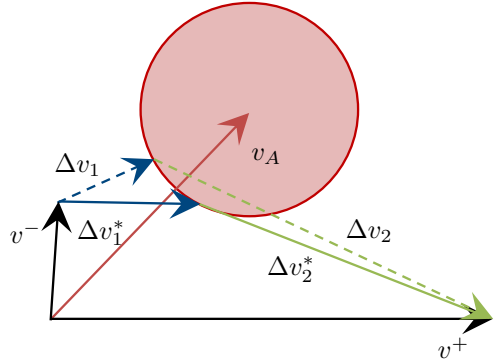

In this step $ \Delta v_{A2A} $ computation is performed in the following simplified way. Let $\boldsymbol{v}^-$ indicate the arrival velocity at an asteroid provided by the Lambert problem. The arrival $\Delta \boldsymbol{v}_1$ is thus the minimum impulse to meet the relative velocity constraint of 2 km/s, as shown in Fig.~\ref{fig:velocity_triangles}. This is achieved by computing the {smallest \(\hat{x}\) for which}  $||\boldsymbol{v}_A - (\boldsymbol{v}^- + x(\boldsymbol{v}_A - \boldsymbol{v}^-)) || \le 2$ km/s. Then $\Delta \boldsymbol{v}_1 =$ ${\hat{x}}(\boldsymbol{v}_A - \boldsymbol{v}^-)$ and $\Delta \boldsymbol{v}_2 = \boldsymbol{v}^+ - (\boldsymbol{v}^-+\Delta \boldsymbol{v}_1)$. Although this approach can lead to sub-optimal $ \Delta v_{1}  + \Delta v_{2} $ (as shown in Fig.~\ref{fig:velocity_triangles} by comparison with the starred quantities), it allows us to compute each $ \Delta v_{A2A}  $ without knowing the subsequent asteroid. {Section~\ref{sec:R} will refine $ \Delta v_{A2A} $, hence completing the description of Fig.~\ref{fig:velocity_triangles}.}

\subsection{Problem Transcription}\label{Sec:CAS}

{A nested-loop optimization \cite{Englander2012a} is employed to find the optimal design parameters $ (\Delta t_{E2A},\Delta t_{A2A}, a_{D})$. The \lq\lq outer loop\rq\rq\ chooses the optimal parameters $ (\Delta t_{E2A},\Delta t_{A2A}, a_{D})$ using a GA optimizer, while the \lq\lq inner loop\rq\rq\ determines the motherships. The BS is used as the inner loop optimizer, building the best campaign of ten mothership for a given set of parameters $ (\Delta t_{E2A},\Delta t_{A2A}, a_{D})$.}


BS is used ten times, one for each mothership. The first mothership is launched at the first available epoch, and the other ones are separated by $ 36.525 $ days, each time removing previously visited asteroids from the main data-set, to prevent repeating objects along the mission. Such an arrangement achieves an even distribution of the initial positions of the ten motherships along the Earth's orbit. The $ BW $ is set to 30, allowing a good compromise between computational effort and solution quality. The quality of each mothership is assessed by the cost function:
\begin{equation} \label{eq:PI_sub}
J_i = \frac{10^{-10} \cdot m_i}{a_{D}^{2} (1+\Delta v_{i}/50)^2}.
\end{equation}
In Eq.~\eqref{eq:PI_sub} $m_i$ denotes the approximate total mass delivered to a target Dyson ring $a_{D}$ (using Edelbaum's formula to estimate the time of flight \cite{edelbaum}) by all the asteroids visited by the $i$-th mother ship (with $ i = 1,...,10 $), and $\Delta v_{i}$ denotes the total maneuver cost of the $i$-th mothership.

After the best ten motherships are generated according to Eq.~\eqref{eq:PI_sub} and for the given choice of $ (\Delta t_{E2A},\Delta t_{A2A}, a_{D}) $, the ultimate overall objective function
\begin{equation} \label{eq:PI}
J = \frac{10^{-10} \cdot m^{tot}}{a_{D}^{2} \sum\limits_{i=1}^{10} (1+\Delta v_{i}/50)^2} 
\end{equation}
is calculated, in which $m^{tot}$ denotes the approximate total mass delivered to the Dyson ring by all the asteroids. The optimization outcome was $\Delta t_{E2A}^* \approx 349$ days, $\Delta t_{A2A}^* \approx 180$ days, $a_{D}^* \approx 1.11$ AU, enabling the transcription of the problem in combinatorial form. 
The optimization was run on the University of Surrey's High Performace Computer (HPC) with 20 cores on a 2.1GHz Intel Xeon processor and took on average 8 hours per generation for a population of 200.  {Algorithm ~\ref{alg:mega-GA} (lines 1-7) provides details of the nested-loop optimization employed to estimate $ (\Delta t_{E2A},\Delta t_{A2A}, a_{D})$.} 

\subsection{Combinatorial Optimization}\label{Sec:CO}

{Once the nested-loop optimization discussed in Sec.~\ref{Sec:CAS} has provided the values for $ (\Delta t_{E2A},\Delta t_{A2A}, a_{D})$, a final BS optimization is run with a much higher beam width $ BW = 600 $. To further enhance the exploration of the combinatorial search space, variations in the $\Delta t_{A2A}$ as estimated by the nested optimization are permitted. As such, Lambert arc transfers connecting two consecutive asteroids are computed 3 times, i.e. for $ \Delta t_{A2A} - 30  \mbox{ days} $, $\Delta t_{A2A}$ and $\Delta t_{A2A} + 30 \mbox{ days}$.} The objective is given by Eq.~\eqref{eq:PI_sub}.

The BS-based combinatorial optimization took on average 2-3 hours to compute 10 mothership chains on a laptop with a 2.8 GHz Quad-Core Intel Core i7 processor (MATLAB environment). The overall approach for the computation of the mothership chains is reported in Algorithm ~\ref{alg:mega-GA}. 

\begin{algorithm*}[!ht]
\caption{Mothership Chains Generation}
\begin{algorithmic}[1]
\setstretch{1}

\State Use GA to optimize the following:

\Function{Fitness\underline{  }ChainBuilding}{$\Delta t_{E2A}$, $\Delta t_{A2A}$, $a_{D}$} 

    \For{each mothership $i = 1$ to $10$}
    
        \State use BS ($ BW = 30 $) minimizing the index $ J_{i} $ \Comment{See Eq.(\ref{eq:PI_sub})}
        
    \EndFor

    \State Calculate fitness value $J$ \Comment{See Eq.(\ref{eq:PI})}

\EndFunction

\State Provided ($\Delta t_{E2A}$, $\Delta t_{A2A}$, $a_{D}$) from GA optimization, apply the relaxation $ \Delta t_{A2A} \pm 30 $ days

\For{each mothership $i = 1$ to $10$}
    
    \State use BS ($ BW = 600 $) minimizing the index $ J_{i} $ \Comment{See Eq.(\ref{eq:PI_sub})}
        
\EndFor

\end{algorithmic}
\label{alg:mega-GA}
\end{algorithm*}

\subsection{Refinement of Encounter Epochs and Flyby Impulses}\label{sec:R}

The last step of the asteroid chain construction is a refinement of the design parameters. This step is necessary to mitigate the effect of the problem transcription \cite{Bellome2021} as described in Sec.~\ref{Sec:CO}. In particular, two post-processing steps are employed:

\begin{itemize}
    \item[1)] Refinement of asteroids encounter epochs and Dyson ring parameters estimate ($\Delta v_{R1}$ in Table~\ref{table:MS_from_BS_test})
\end{itemize}

Once the asteroids sequences are known for each mothership, the performance index $ J_{i} $ from Eq.~\eqref{eq:PI_sub} becomes a function of only continuous-varying parameters
\begin{equation} \label{eq:refinement}
    J_{i} = f(\hat{S}, T, a_{D}, i_{D}, \Omega_{D}),
\end{equation}
where $\hat{S} $ are the known asteroids sequence, and vector $ T $ encodes the encounter epochs for each object. 

A PSO was first employed {to find T} that minimizes $ J_{i}$ for each mothership. This optimization was followed by one in which the objective function \eqref{eq:PI} was maximized acting only on the ring parameters $(a_{D}, i_{D}, \Omega_{D})$. In this process, the arrival asteroids masses were still estimated with Edelbaum's approach \cite{edelbaum}. 
\begin{itemize}
    \item[2)] Refinement of flyby impulses ($\Delta v_{R2}$ in Table~\ref{table:MS_from_BS_test})
\end{itemize}

Suppose the arrival $\boldsymbol{v}^-$ and final velocity $\boldsymbol{v}^+$ at an asteroid were available from Lambert's arcs, the flyby impulse optimization consists of finding the $\Delta \boldsymbol{v}_1$ and $\Delta \boldsymbol{v}_2$ that 
minimize
\begin{equation}
    J = \Delta {v}_1 + \Delta {v}_2 
\end{equation}
subject to
\begin{align}
    ||\boldsymbol{v}^-+ \Delta \boldsymbol{v}_1 - \boldsymbol{v}_A|| & \le  2\  \textrm{km/s} \\
    \boldsymbol{v}^- + \Delta \boldsymbol{v}_1 + \Delta \boldsymbol{v}_2 & =  \boldsymbol{v}^+.
\end{align} 
As illustrated in Fig.~\ref{fig:velocity_triangles}, the optimal values $\Delta \boldsymbol{v}^*_1$ and $\Delta \boldsymbol{v}^*_2$ can produce significant improvement in the total $\Delta v$ compared to $\Delta \boldsymbol{v}_1$ that minimizes the impulse satisfying $ ||\boldsymbol{v}^-+ \Delta \boldsymbol{v}_1 - \boldsymbol{v}_A|| \le  2\  \textrm{km/s}$. 
Table~\ref{table:MS_from_BS_test} reports the incremental savings in total $\Delta v$ achieved by the two refinements. These improvements required a few minutes to run on a laptop with a 2.8 GHz Quad-Core Intel Core i7 processor (MATLAB environment).

\section{Rendezvous Table Generation}\label{sec:TG} 
In this section, the optimal flight time for a given asteroid to any of the 12 Dyson ring stations is estimated. First, a time-optimal trajectory is obtained from each asteroid's original orbit to the Dyson ring orbit (no rendezvous).  Then, phasing-based interpolation is used to estimate the time taken for the said asteroid to reach a given station. 

\subsection{Single Asteroid Transfer to the Dyson Orbit}\label{SAT_m}
The flight time of each asteroid must be minimized, as an increase in flight time corresponds to a decrease in the asteroid's mass upon arrival. 
To this end, optimal control theory and indirect optimization methods are utilized to obtain optimal flight times. \par 
A relatively accurate estimate of the initial costates is required when solving the time-optimal problem through indirect methods. This is obtained by first solving the optimal energy problem for the minimum time of flight calculated using Edelbaum's equation \cite{edelbaum}.

\subsubsection{Energy-Optimal Problem}
The energy-optimal nonlinear control problem can be formulated as: minimize
\begin{equation}
    J = \frac{1}{2}f_{\text{ATD}} \int^{t_f}_{t_0} \tau^2 \ dt 
\end{equation}
subject to 
\begin{equation} 
    \dot{\boldsymbol{x}} = \boldsymbol{A}(\boldsymbol{x})+ \boldsymbol{B}(\boldsymbol{x})f_{\text{ATD}}\tau \hat{\boldsymbol{\alpha}},
    \label{stateeq}
\end{equation}
where $\tau$ and $\hat{\boldsymbol{\alpha}}$ represent the engine unbounded thrust ratio and thrust direction unit vector, respectively. $f_{ATD}$ indicates the fixed magnitude of acceleration provided by the ATD. The modified equinoctial elements are used in the problem formulation, thus the expressions for $ \boldsymbol{A}(\boldsymbol{x})$ and $\boldsymbol{B}(\boldsymbol{x})$ are those given in \cite{A&B}, and $\boldsymbol{x} = [p, f,g,h,k,L]$.  \par 
To solve the nonlinear energy-optimal problem using indirect methods,
a costate vector $\boldsymbol{\lambda} \in \mathbb{R}^6$ is introduced. The Hamiltonian is 
\begin{equation}
    \mathcal{H} =\frac{1}{2} f_{\text{ATD}} \tau^{2}+\boldsymbol{\lambda} ^{T} \boldsymbol{A}(\boldsymbol{x})+\boldsymbol{\lambda} ^{T}  \boldsymbol{B}(\boldsymbol{x}) f_{\text{ATD}} \tau \widehat{\boldsymbol{\alpha}}
\end{equation}
and the optimal control is derived by minimizing $\mathcal{H}$:
\begin{equation}
    \widehat{\boldsymbol{\alpha}}^{\star}=-\frac{\boldsymbol{B}^{T} \boldsymbol{\lambda}}{\left\|\boldsymbol{B}^{T} \boldsymbol{\lambda}\right\|}, \quad \tau^{\star}=\left\|\boldsymbol{B}^{T} \boldsymbol{\lambda}\right\|.
    \label{control}
\end{equation}
Substituting Eq.~\eqref{control} into Eq.~\eqref{stateeq}, the equations of motion and the Euler-Lagrange equation with control can be obtained as follows: {
\begin{equation} 
\begin{gathered}
\dot{\boldsymbol{x}}=\boldsymbol{A}(\boldsymbol{x}) + \boldsymbol{B}(\boldsymbol{x})   f_{\text{ATD}}  \tau^{\star}\widehat{\boldsymbol{\alpha}}^{\star}  \\
\dot{\boldsymbol{\lambda}}=-\frac{\partial \left(\boldsymbol{\lambda}^{T} \boldsymbol{A}(\boldsymbol{x})\right)}{\partial \boldsymbol{x}} - f_{\text {ATD}} \tau^{\star}\  \frac{\partial \left( \boldsymbol{\lambda} ^{T} \boldsymbol{B}(\boldsymbol{x})\widehat{\boldsymbol{\alpha}}^{\star}\right)}{\partial \boldsymbol{x}} .\end{gathered}
\label{el}
\end{equation}
For simplicity the second term in the equation of the Lagrangian multipliers is neglected in the energy-optimal problem as we are only using this solution to obtain a first guess solution. The full expression of the Euler-Lagrange equations is used later for the time-optimal problem.}

{ }

This problem can now be solved using the shooting method. Since the true longitude is treated as a free variable, the shooting function is: 
\begin{equation}
   \boldsymbol{\phi}\left(\boldsymbol{\lambda}\left(t_{0}\right)\right)=\left[\boldsymbol{x}_{1:5}^{T}\left(t_{f}\right)-{\boldsymbol{x}_{f}}_{1:5}^{T}, {\lambda}_L(t_f)\right]^{T}=\boldsymbol{0}.
   \label{shooting1}
\end{equation}
As mentioned, a non-intuitive initial guess for the initial costates is required for the shooting method. This is procured by linearizing Eq.~\eqref{el} and obtaining an analytical solution for the initial costates of the linearized system using the method described in \cite{diwu}. These linearized initial costates are then used as an initial guess to solve the nonlinear energy-optimal problem. {To effectively use the energy-optimal solution as a first guess for the time-optimal problem, the time of flight $\Delta t$ in the energy-optimal is adjusted such that obtained $\Delta v$ is equivalent to $ f_{\text{ATD}}\Delta t$.}

\subsubsection{Time Optimal Problem}\label{toptfree}
This section describes the method used to obtain the time optimal trajectory from the original asteroid orbit to the Dyson ring orbit. This problem consists in minimizing 
\begin{equation}
    J =\int^{t_f}_{t_0} dt, 
\end{equation}
subject to
\begin{equation} 
    \dot{\boldsymbol{x}} = \boldsymbol{A}(\boldsymbol{x})+ \boldsymbol{B}(\boldsymbol{x}) f_{\text{ATD}}\tau  \widehat{\boldsymbol{\alpha}}.
    \label{stateeq1}
\end{equation}
Note that since the problem assumes constant thrust acceleration, the thrust ratio is assumed to be unity. A costate vector is introduced and the Hamiltonian 
\begin{equation}\label{hamil}
    \mathcal{H} =1+\boldsymbol{\lambda} ^{T} \boldsymbol{A}(\boldsymbol{x})+\boldsymbol{\lambda} ^{T} \boldsymbol{B}(\boldsymbol{x}) f_{\text{ATD}}\tau \widehat{\boldsymbol{\alpha}}
\end{equation}
is established.
{The thrust direction that minimizes $\mathcal{H}$ is the same as for the energy-optimal problem, while $\tau^{\star} = 1$ at all times.  The Euler-Lagrange equations are given by Eq. \eqref{el}.
The problem can be solved using a shooting method with the initial costates of the nonlinear energy-optimal solution as initial guesses.} The shooting function in this case is defined as follows: 
\begin{equation}
    \phi\left(\boldsymbol{\lambda}_{t_{0}}, t_{f}\right)=\left[\boldsymbol{x}_{1:5}^{T}\left(t_{f}\right)-{\boldsymbol{x}_{f}}_{1:5}^{T}, \boldsymbol{\lambda}_L(t_f), H\left(t_{f}\right)\right]=0,
\end{equation}
{The final constraint on the Hamiltonian derives from 
\begin{equation}
H(t_f) -\boldsymbol{\lambda}_{1:5}(t_f)^T \dot{\boldsymbol{x}}_{f_{1:5}} = 0,    
\end{equation}
in which $\dot{\boldsymbol{x}}_{f_{1:5}} =\boldsymbol{0}$ as the five targeted orbital elements are constant}.

\subsection{Estimating the Time of Flight for Each Station}\label{sec:ETF}
This section aims to determine all possible trajectories between a single asteroid to any Dyson ring station. The procedure is divided into three steps for clarity: 
\begin{enumerate}
    \item Primary data collation: calculating the time-optimal trajectories within one asteroid orbital period following the flyby.
    \item Continuation: obtaining time-optimal trajectories with starting epochs beyond one asteroid orbital period following the flyby.
    \item Phase matching: interpolating phase differences to obtain estimated trajectories for each station. 
\end{enumerate}

\subsubsection{Primary Data Collation}
The procedure for the primary data collection can be summarized as:
\begin{enumerate} \label{listref_primary}
\item  Calculate the orbital period (P) of the asteroid, and obtain {
$t_{0,primary} = t_{flyby} + 30 d + (k-1)(\frac{P}{n})$
where, $k = 1:n$.  During the competition $n =16$ was used}, as calculating the trajectories for 16 points within each asteroid orbital period proved to be a good trade-off between computational time and solution accuracy. (Note that 30 days are added to the time of the mothership flyby to account for the time required for the ATD activation.)  

\item For each start time in $t_{0,primary}$:
\begin{enumerate}
    \item Solve the time-optimal problem described in Sec.~\ref{SAT_m}. 
    \begin{itemize}
        \item For $k =1$, obtain an initial guess of the costates for the time optimal problem by solving the energy-optimal problem. 
        \item For $k = 2:n$, use the initial costates from the $(k-1)$ time step as the initial guess for the $(k)$th time step.
    \end{itemize}
    \item Store the arrival time $t_{f,primary}$, and the initial costates $\boldsymbol{\lambda}_{t_0}$. (Note that $\Delta t _{primary} = t_{f,primary} - t_{0,primary} $ )
     \item \label{deltaLcalc} Repeat the following for station $= 1:12$ 
    \begin{enumerate}[i)]
    \item Calculate the station's true longitude ($L_{S}$) at $t_{f,primary}$. \footnote{ {Note that the phase of the first building station was assumed to be zero at the starting epoch of 95739 MJD. No optimization was conducted on this variable throughout the competition.}}
    \item Obtain the asteroid's true longitude ($L_{A}$) at $t_{f,primary}$.
    \item Calculate the phase difference upon arrival: $\Delta L = L_{A} -  L_{S}$   
    \item Store $\Delta L$ for each station. 
    \end{enumerate}
\end{enumerate}
 
\end{enumerate}

\begin{table}[H]
\caption{Phasing table data for primary epochs}

\centering 
\begin{tabular}{cccccc}
\hline
\multirow{2}{*}{\begin{tabular}[c]{@{}c@{}}Departure\\ time\end{tabular}} & \multicolumn{3}{c}{Phase difference}                                                                    & \multirow{2}{*}{\begin{tabular}[c]{@{}c@{}}Arrival \\ time\end{tabular}} & \multirow{2}{*}{\begin{tabular}[c]{@{}c@{}}Initial \\ costates\end{tabular}} \\ \cline{2-4}
                                                                          & \multicolumn{1}{c|}{1}                & \multicolumn{1}{c|}{…} & \multicolumn{1}{c}{12}                &                                                                          &                                                                              \\ \hline
$t_{0,primary} (1)$                                                       & \multicolumn{1}{c|}{$\Delta L(1,1)$}  & \multicolumn{1}{c|}{…} & \multicolumn{1}{c|}{$\Delta L(1,12)$}  & $t_{f,primary}(1)$                                                       & $\boldsymbol{\lambda}_{t_0}(1) $                                           \\
$t_{0,primary} (2) $                                                      & \multicolumn{1}{c|}{$\Delta L(2,1)$}  & \multicolumn{1}{c|}{…} & \multicolumn{1}{c|}{$\Delta L(2,12) $} & $t_{f,primary}(2)$                                                       & $\boldsymbol{\lambda}_{t_0}(2)$                                            \\
…                                                                         & \multicolumn{1}{c|}{..}               & \multicolumn{1}{c|}{…} & \multicolumn{1}{c|}{…}                 & …                                                                        & …                                                                            \\
$t_{0,primary}(16)$                                                       & \multicolumn{1}{c|}{$\Delta L(16,1)$} & \multicolumn{1}{c|}{…} & \multicolumn{1}{c|}{$\Delta L(16,12)$} & $t_{f,primary} (16)$                                                     & $\boldsymbol{\lambda}_{t_0 }(16) $                                         \\ \hline
\end{tabular}
\label{pt1}
\end{table}

 Table~\ref{pt1} shows the output of the primary data collection procedure. {The information of phase difference is collated as it informs us about the distance between any given station and the asteroid upon arrival. This information can then be used to determine the adjustments that need to be made in the departure time for a particular asteroid to rendezvous with a particular station. }

\subsubsection{Continuation}
{If an asteroid is visited in the early phases of the mission, there are multiple transfer opportunities to the stations beyond its first orbital period following the flyby. The following steps are taken to calculate these opportunities. Note that this steps does not require computing any new trajectories as for time-optimal orbit transfer the geometry of the problem repeats after one orbital period. } 

 \begin{enumerate}
     \item While the $\max{(t_f)}$ is within the mission timeline:
     \begin{enumerate}
         \item Define a secondary set of departure times:\\  $t_{0,secondary} = t_{0,primary} + P$
          \item For each start time in $t_{0,secondary}$:
         \begin{enumerate}[i)]
             \item As the problem geometry for the secondary set remains the same as the primary one \\
          $ \Delta t_{secondary} = \Delta t_{primary} \ \text{and} \ \boldsymbol{\lambda}_{secondary} = \boldsymbol{\lambda}_{primary}$
          \item Calculate the arrival time:\\ 
          $ t_{f, secondary} =   t_{0,secondary} + \Delta t_{secondary}$
          \item Calculate the phase difference upon arrival at each station, $\Delta L$.
         \end{enumerate}
         
     \end{enumerate}
 \end{enumerate}
 {Once all free fast variable trajectories for a given asteroid are calculated, the phasing table can be completed, as shown in Table~\ref{pt2}. Note that while the initial costates for primary and secondary sets repeat identically, the node differences $\Delta L$ do not repeat as the asteroid and the Dyson ring do not share the same orbital period. The values of $\Delta L$ are the key data to compute the rendezvous trajectories, as shown in the next section.}
 
\begin{table}[H]
\caption{Phasing table data for primary and secondary epochs}
\centering 
\begin{tabular}{cccccc}
\hline
\multirow{2}{*}{\begin{tabular}[c]{@{}c@{}}Departure\\ time\end{tabular}} & \multicolumn{3}{c}{Phase difference}                                                & \multirow{2}{*}{\begin{tabular}[c]{@{}c@{}}Arrival \\ time\end{tabular}} & \multirow{2}{*}{\begin{tabular}[c]{@{}c@{}}Initial \\ costates\end{tabular}} \\ \cline{2-4}
                                                                          & \multicolumn{1}{c|}{1}                & \multicolumn{1}{c|}{…} & 12                 &                                                                          &                                                                              \\ \hline
$t_{0,primary} (1)$                                                       & \multicolumn{1}{c|}{$\Delta L(1,1)$}  & \multicolumn{1}{c|}{…} & $\Delta L(1,12)$   & $t_{f,primary}(1)$                                                       & $\boldsymbol{\lambda}_{t_0}(1) $                                           \\
$t_{0,primary} (2) $                                                      & \multicolumn{1}{c|}{$\Delta L(2,1)$}  & \multicolumn{1}{c|}{…} & $\Delta L(2,12) $  & $t_{f,primary}(2)$                                                       & $\boldsymbol{\lambda}_{t_0}(2)$                                            \\
…                                                                         & \multicolumn{1}{c|}{..}               & \multicolumn{1}{c|}{…} & …                  & …                                                                        & …                                                                            \\
$t_{0,primary}(16)$                                                       & \multicolumn{1}{c|}{$\Delta L(16,1)$} & \multicolumn{1}{c|}{…} & $\Delta L(16,12)$  & $t_{f,primary} (16)$                                                     & $\boldsymbol{\lambda}_{t_0 }(16) $                                         \\ \hline
$t_{0,secondary} (1)$                                                     & \multicolumn{1}{c|}{$\Delta L(17,1)$} & \multicolumn{1}{c|}{…} & $\Delta L(17,12)$  & $t_{f,secondary}(1)$                                                     & $\boldsymbol{\lambda}_{t_0}(1) $                                           \\
$t_{0,secondary} (2) $                                                    & \multicolumn{1}{c|}{$\Delta L(18,1)$} & \multicolumn{1}{c|}{…} & $\Delta L(18,12) $ & $t_{f,secondary}(2)$                                                     & $\boldsymbol{\lambda}_{t_0}(2)$                                            \\
…                                                                         & \multicolumn{1}{c|}{..}               & \multicolumn{1}{c|}{…} & …                  & …                                                                        & …                                                                            \\
$t_{0,secondary}(16)$                                                     & \multicolumn{1}{c|}{$\Delta L(32,1)$} & \multicolumn{1}{c|}{…} & $\Delta L(32,12)$  & $t_{f,secondary} (16)$                                                   & $\boldsymbol{\lambda_{t_0 }}(16) $                                         \\ \hline
…                                                                         & \multicolumn{1}{c|}{..}               & \multicolumn{1}{c|}{…} & …                  & …                                                                        & …                                                                            \\ \hline
\end{tabular}
\label{pt2}
\end{table}

\subsubsection{Phase Matching}
{
To rendezvous with a station, the true longitude of the asteroid and the station must coincide. All the rendezvous opportunities of one asteroid with all the stations are accurately estimated by interpolating the data in }Table~\ref{pt2}. Cubic spline interpolation is used to obtain the departure time, arrival time and initial costates for which $\Delta L = 2k \pi, k \in \mathbb{Z}$. \par
\begin{figure}[H]
    \centering
    \begin{minipage}{.48\textwidth}
        \centering
        \includegraphics[width = \textwidth]{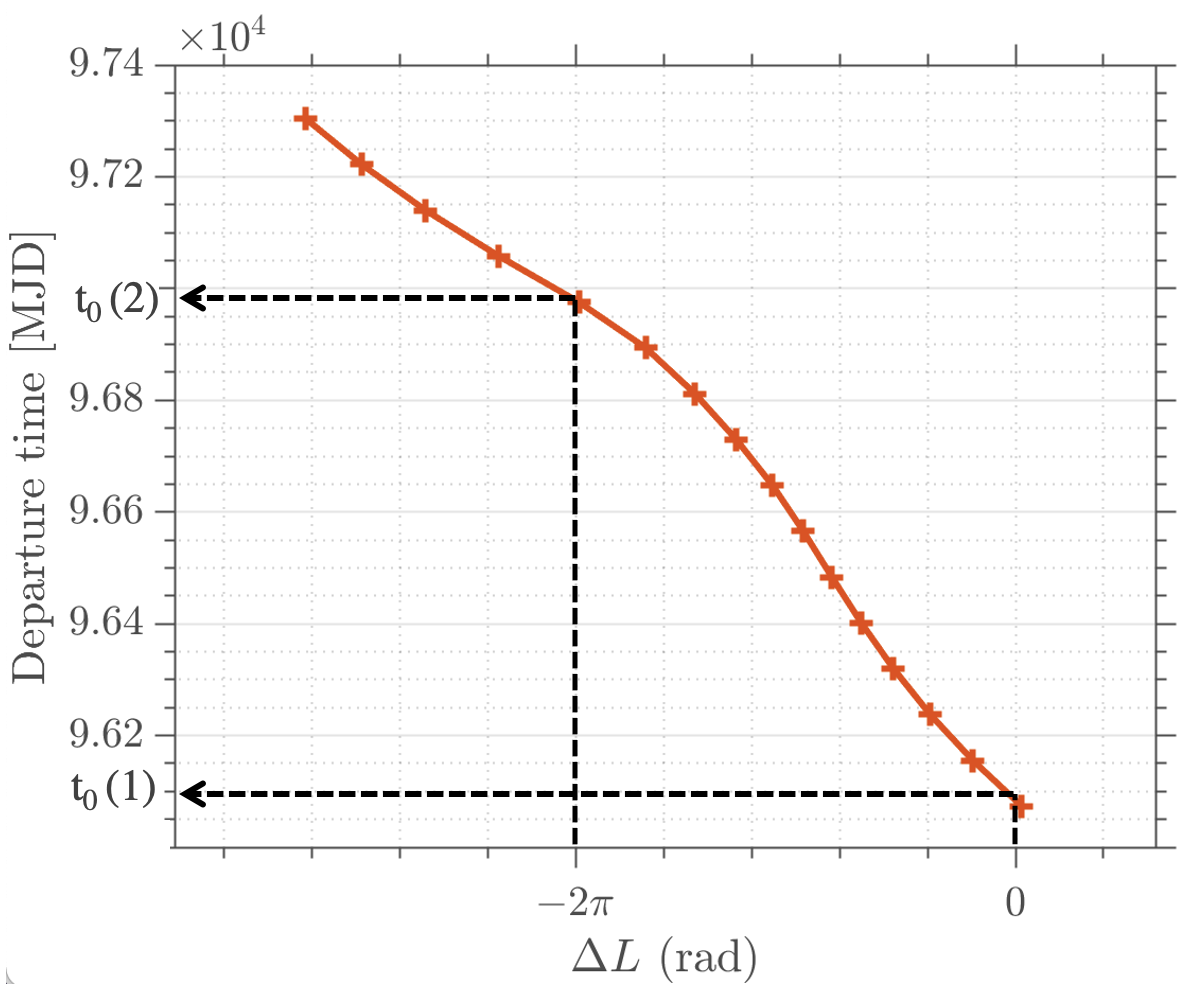}
        \caption{Determination of two departure times (within the first orbit period) for a rendezvous trajectory with Station 1}
        \label{fig:depature}
    \end{minipage}%
    \quad 
    \begin{minipage}{0.48\textwidth}
        \centering
         \includegraphics[width = \textwidth]{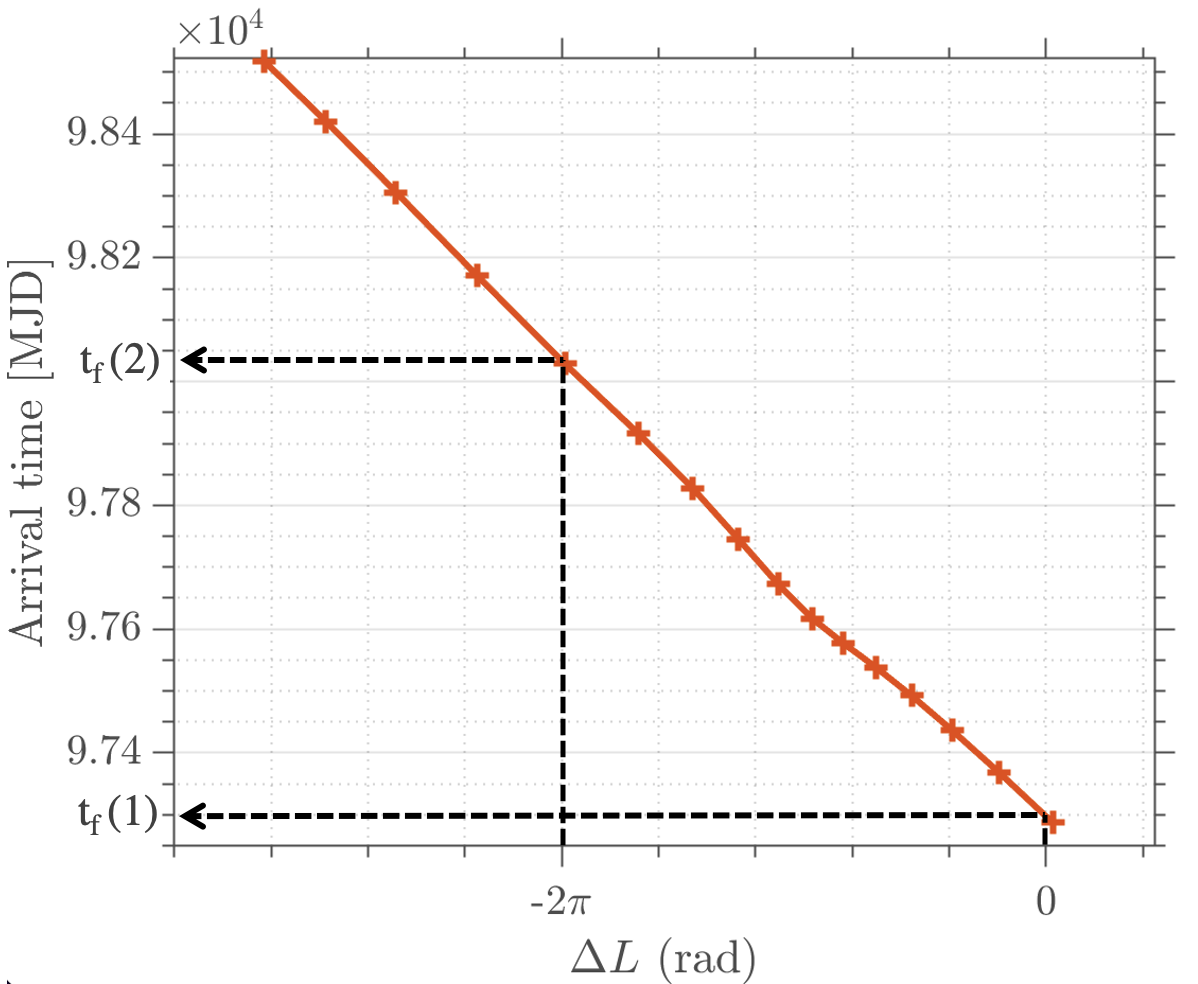}
        \caption{Determination of two arrival times for a rendezvous trajectory with Station 1, for departure timed from Fig. \ref{fig:depature}}
        \label{fig:arrival}
    \end{minipage}
\end{figure}

 
Figures \ref{fig:depature} and \ref{fig:arrival} {show the interpolation process to obtain the departure and arrival times for an asteroid to get to Station 1, using the primary data (i.e., within the first orbital period after flyby). The red lines are obtained through cubic spline interpolation of the 16 points included in the primary data. Two rendezvous opportunities are identified when $\Delta L = 0$ and $\Delta L  = -2\pi$. }

Following the interpolation process, the mass upon arrival at each station is calculated from $ m_f=m_{0}^{ast}-\alpha m_{0}^{ast} (\Delta t)$. 

{
For a single asteroid, this process is repeated for all allowed orbital periods (as in Table \ref{pt2}) and for all 12 stations. The outcome is stored in the struct given in Table~\ref{m1} that contains a sequence of possible starting epochs, time of flights, mass at arrival, and initial costates obtained from interpolation. By repeating this process for all asteroids in a given mothership chain, the rendezvous table (Table~\ref{m2}) is generated. This table is then taken as input for the dispatcher to determine the asteroid trajectories that optimize the performance index. }

\begin{table}[H]
          \caption{Data in a single structure of the table (i.e. all transfers from a single asteroid to a single station). {Note that ellipsis are used to denote the existence of additional entries in the sequences.}}
          \label{m1}
\centering 
\begin{tabular}{cc}
\hline
Struct$(i,j)$   & $A_i2S_j $ \\\hline
$t_0$          & $t_0(1), t_0(2), \dots$ \\
$t_f$         & $t_f(1), t_f(2), \dots$\\
$m_f$          & $m_f(1), m_f(2), \dots$\\
$\boldsymbol{\lambda}$      &$\boldsymbol{\lambda}_{t_0}(1), \boldsymbol{\lambda}_{t_0}(2), \dots$ \\
\hline
        \end{tabular}

\end{table}

\begin{table}[H]
         \caption{Overview of the rendezvous table, taken as an input from the dispatcher. {Note that ellipsis are used to denote that the table of structs extends over all asteroids and all 12 stations.}}
          \label{m2}
\centering 
        \begin{tabular}{|c|c|c|c}
\hline
   & Station 1 & Station 2            & ... \\
\hline
Asteroid 1 &  $A_12S_1$ & $A_12S_2$ &  ...   \\
Asteroid 2 &  $A_22S_1$ &  $A_22S_2$       & ...    \\
...        &  ... &       ...      &   ...\\ 
        \end{tabular}
\end{table}

To reduce the computational time, the algorithm for computing the time-optimal transfers was coded in C++ and solved with simple shooting using the multidimensional root-finding algorithm implemented in the GNU scientific library \cite{gough2009gnu}. The required derivatives were computed automatically using the DACE library \cite{rasotto2016differential}. The table generation for the optimal solution required about 12 minutes on a laptop with a 1.9 GHz Quad-Core Intel Core i7 processor. 

\section{The Dispatcher} \label{sec:4.1}
The objective of the dispatcher is to allocate asteroids to stations optimally. As an input, the algorithm takes the table described in Sec.~\ref{sec:TG} containing the list of available asteroids and all their possible arrival times and arrival masses for a given ring altitude (see Table~\ref{m2}). The approach is then to find the sequence of stations and asteroids used to build each station that maximizes the GTOC objective function, Eq.~\eqref{eq:obj}, while satisfying the constraint on the minimum time between the latest arrival at a station and the earliest at the following station (constraint (\ref{constraint90days}) in Sec.~\ref{introduction}). Since the input data fixes the ring radius, $a_{\text{D}}$, two terms remain in Eq.~\eqref{eq:obj}. Concerning the minimum mass, $m_{\mathrm{min}}^{\text{S}}$, the optimal condition is achieved by allocating as many asteroids as possible such that the total mass in each station $m^{S_i}$, where $i = 1,\ldots,12$, is distributed evenly. The $\Delta v_{i}$ term is reduced by trimming any mothership chains whenever an asteroid at the end of the chain is not allocated. 

In order to do so, an optimization problem was set up. The optimization variables are: $\bm x_S$, 12 real numbers between 0 and 1 for defining the sequence of stations; $\bm x_{NA}$, 12 integers {between 12 and 36 (other boundaries were tested during the competition but these were the best-performing ones)} defining the number of asteroids in each station; $x_{\Delta t}$, a scalar that defines the minimum time interval between two contiguous asteroid arrivals when assigned to two consecutive stations. This time interval was initially fixed during the optimization and set equal to 91 days to allow for a margin with respect to the constraint of 90 days. However, a sensitivity analysis showed significant variability in the optimization results even with minor variations in this parameter. Thus, this parameter was included among the optimization variables.

The objective function's core is divided into three main parts described below. Intermediate results for the optimal solution are shown for illustrative purposes, but the same approach is followed in each function call during the optimization process. 

\subsection{Asteroids' First Allocation} The sequence of stations, $S_i$ with $i =1,\ldots,12$, is defined by ordering the values in $\bm x_S$ and taking the indexes corresponding to the ordered elements. At this point, starting from the first station, $S_1$, the asteroids are ordered based on arrival, with the earliest arrival first. The first $\bm x_{NA}(1)$ asteroids are assigned to $S_1$ and removed from the list. A similar procedure is followed for the second station and the following ones. In this case, however, only arrivals after the latest arrival in the previous station plus $x_{\Delta t}$ are considered. 

\begin{figure}
\begin{center}
\includegraphics[width=0.45\columnwidth]{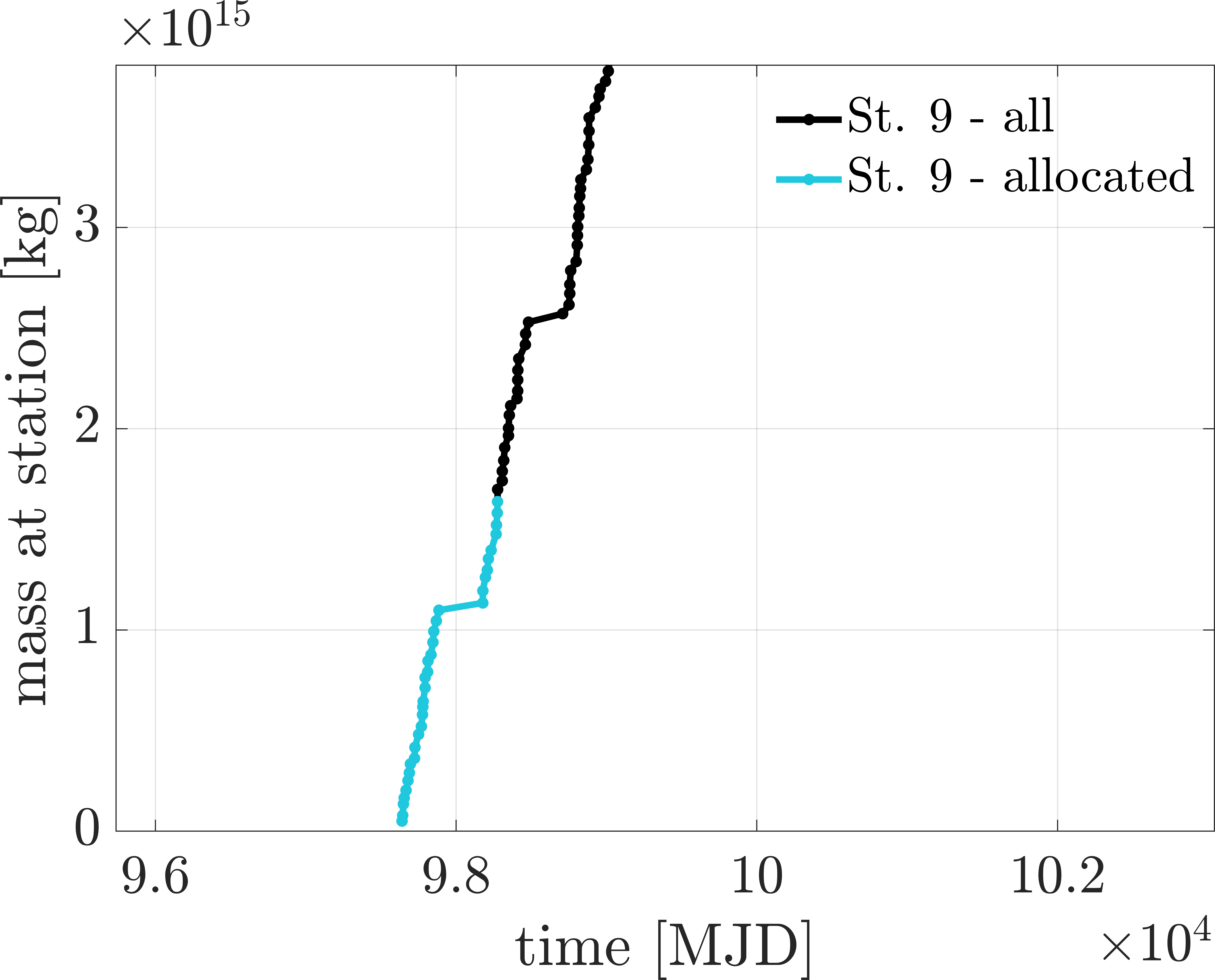}
\includegraphics[width=0.45\columnwidth]{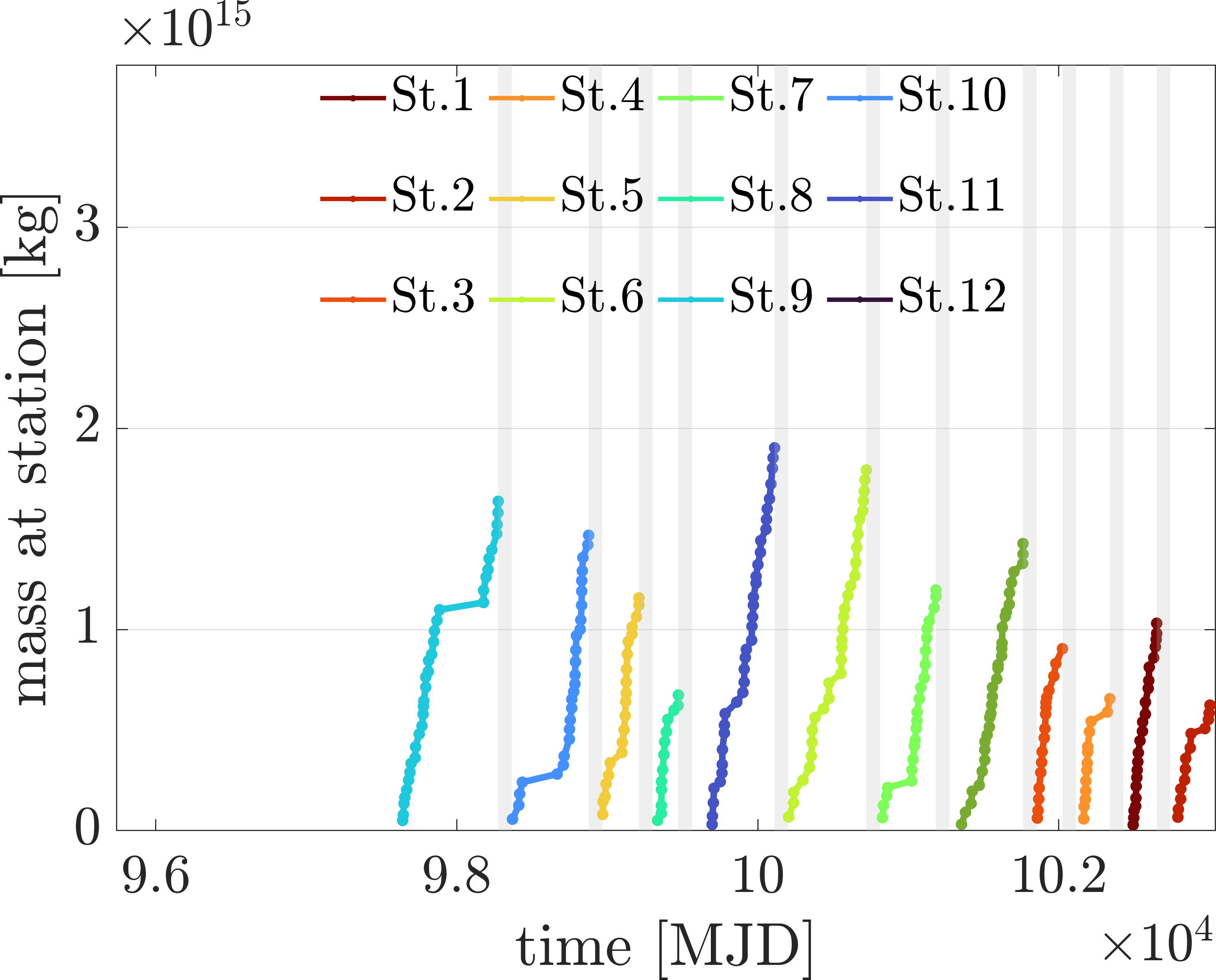}
\caption{{Mass at station as a function of the asteroid arrival time for a single station (left) and for the optimal solution as a result of first asteroid allocation (right).}}
\label{fig:DispatcherFirstAllocation}
\end{center}
\end{figure}

In Fig.~\ref{fig:DispatcherFirstAllocation} (left panel), the mass at the first station, $m^{S_9}$, is shown as a function of the arrival time. The light blue points highlight the asteroids assigned to this station, while the remaining ones are gray. Fig.~\ref{fig:DispatcherFirstAllocation} (right panel) shows the final result, where the gray vertical shades indicate the time interval between two consecutive stations, $x_{\Delta t} =$ 90.34 days. Here, 282 asteroids of the 313 available were assigned to a station, with a minimum mass $m_{\mathrm{min}}^{\text{S}} = 6.2 \times 10^{14}$ kg corresponding to $m^{S_2}$.

\subsection{Asteroids' Distribution Refinement} At this stage the total mass at each station might vary significantly. For instance, in the optimal solution, the difference between the minimum and maximum values (maximum mass corresponds to $m^{S_{11}} = 1.9\times10^{15}$ kg) is greater than $1.2\times10^{15}$ kg. Since the objective function only depends on the minimum value, in principle, all other asteroids sequences could be cut at values equal to or greater than $m_{\mathrm{min}}^{\text{S}}$ without affecting the mass contribution in the objective function. Ideally, the optimal solution would involve all asteroids allocated at the earliest arrival time (i.e. with maximum mass $m_f$, see Table~\ref{m1}), and equally distributed across the stations. That is, $m^{S_i} = \sum_{j} m_{f,Ast_j}/12$, in which $m_{f,Ast_j}$ is the final mass of the $j$th asteroid. Therefore, an iterative strategy is followed, where asteroids are removed from high mass stations and placed in the stations with the lowest mass. This terminates when no more asteroids can be removed or allocated to a station. 

Some key design choices for implementing the algorithm are summarized below. While other options were explored during the competition, these were those providing the best results. 
\begin{enumerate}
\item The algorithm is fully deterministic with no direct dependencies on the input parameters. That is, the starting point is the asteroid distribution resulting from the first asteroid allocation.
\item The average mass across all stations is used as a mass threshold for removing asteroids. An exception is made in the last iteration, where the minimum mass is used instead.
\item Asteroids with low masses are selected first for removal. They are pooled together with the asteroids not assigned during the first asteroid allocation.
\item Only one asteroid, the one characterized by the earliest arrival compatible with the problem constraints, is allocated to the station with minimum mass in each iteration. 
\end{enumerate}

\begin{figure}
\begin{center}
\includegraphics[width=0.45\columnwidth]{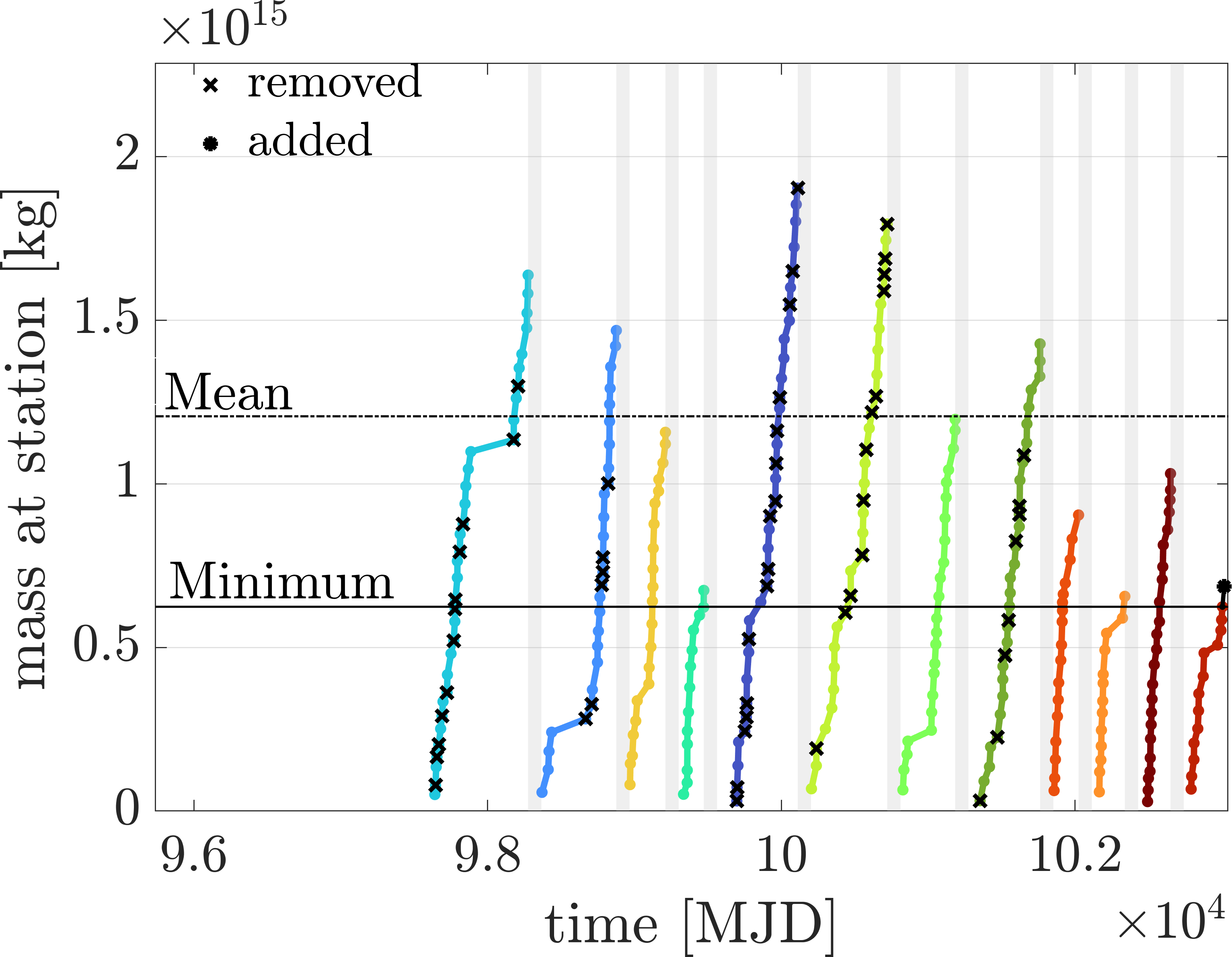}
\includegraphics[width=0.45\columnwidth]{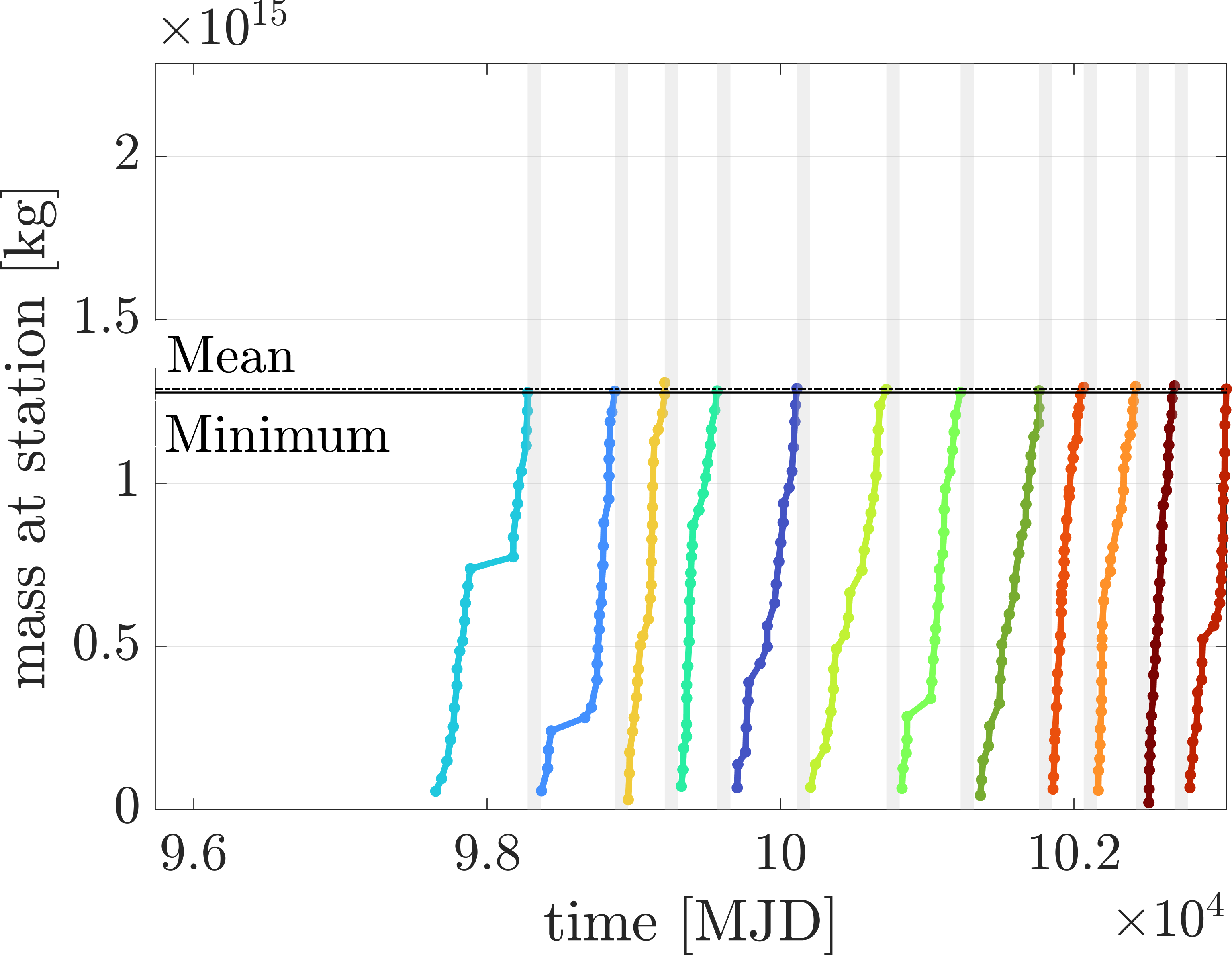}
\caption{Mass at the stations as a function of the asteroid arrival time for the optimal solution as a result of asteroids' distribution refinement.}
\label{fig:DispatcherRefinement}
\end{center}
\end{figure}

\begin{figure}
\begin{center}
\includegraphics[width=0.45\columnwidth]{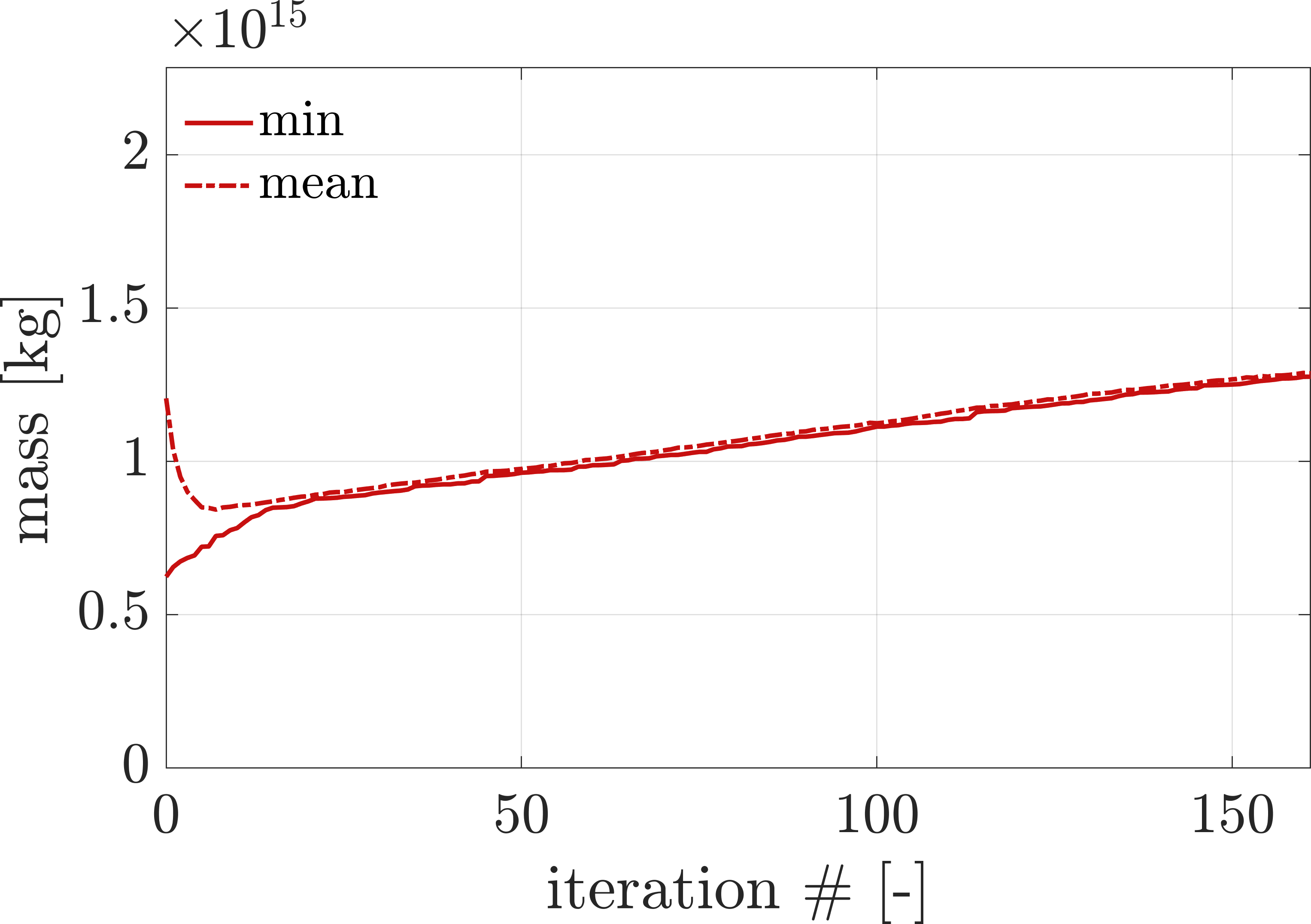}
\includegraphics[width=0.45\columnwidth]{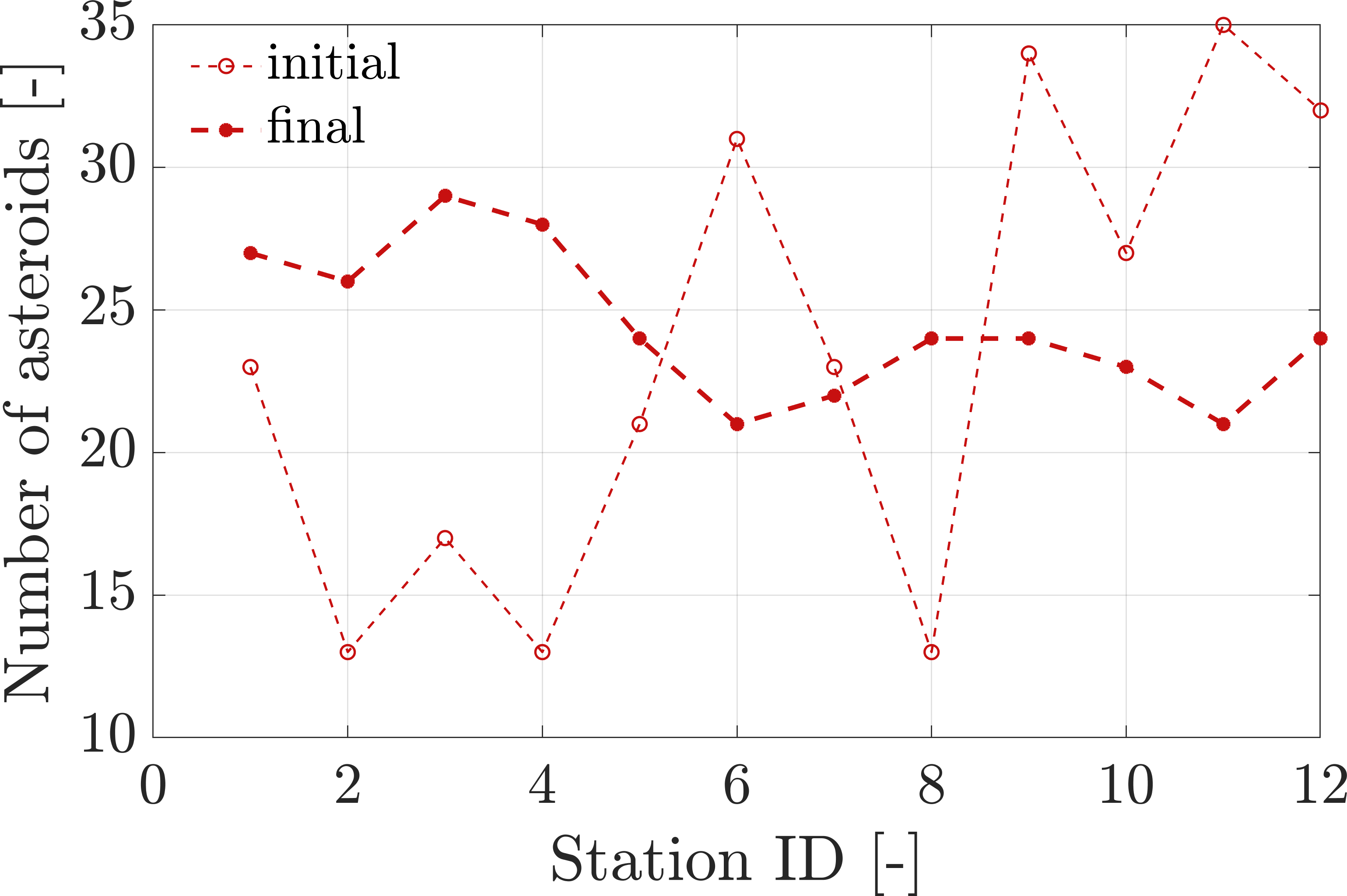}
\caption{Minimum and average mass in all the stations as a function of the iteration of distribution refinement (left panel) and number of asteroids in each station (right panel) at the beginning and at the end.}
\label{fig:DispatcherRefinementConv}
\end{center}
\end{figure}

In Fig.~\ref{fig:DispatcherRefinement} (left panel), the initial mass distribution is shown for the optimal solution, together with the minimum and average mass values (black lines). Here, the black crosses highlight the asteroids removed in the first iteration, and the black asterisk shows the asteroid allocated to $S_2$ in the same iteration. The minimum and average mass across the stations as a function of the iteration is shown in Fig.~\ref{fig:DispatcherRefinementConv} (left panel). The distance between the two curves progressively decreases until the algorithm stops at the 161st iteration, with an average mass just slightly above the minimum one (the difference is around $1.1\times10^{13}$ kg). The final distribution is shown in Fig.~\ref{fig:DispatcherRefinement} (right panel), showing an almost uniform mass distribution across all the stations, resulting in $m_{\mathrm{min}}^{\text{S}} = 1.3 \times 10^{15}$ kg corresponding to $m^{S_9}$. The number of asteroids in each station before and after this refinement step is shown in Fig.~\ref{fig:DispatcherRefinementConv} (right panel). The maximum variability in the number of asteroids per stations reduces from 22 to 8. Importantly, after the distribution refinement, 293 asteroids are assigned, 11 more than at the end of the first asteroid allocation.

\subsection{MS Chain Trimming} Now that all the asteroids are optimally assigned to the stations, a reduction in the motherships $\Delta v$ is achieved by trimming the chains whenever an asteroid not assigned is at the end of one chain. In principle, a further reduction could be obtained if deep space maneuvers (DSMs) were added to skip asteroids. However, this was not implemented within this optimization process but only as a final refinement of the solution (see Sec.~\ref{ssec:DSM}).

In the optimal solution, 19 out of the 20 asteroids not assigned are at the end of the mothership chains. The mothership chain trimming allowed for a reduction of 8.48 km/s in $\Delta v$. Details of the $\Delta v$ after the trimming, indicated with $\Delta v_{trim}$, are reported in Table~\ref{table:MS_from_BS_test}.

\subsection{Computing the Optimal Solutions}
\begin{figure}[H]
\begin{center}
\includegraphics[width=0.5\columnwidth]{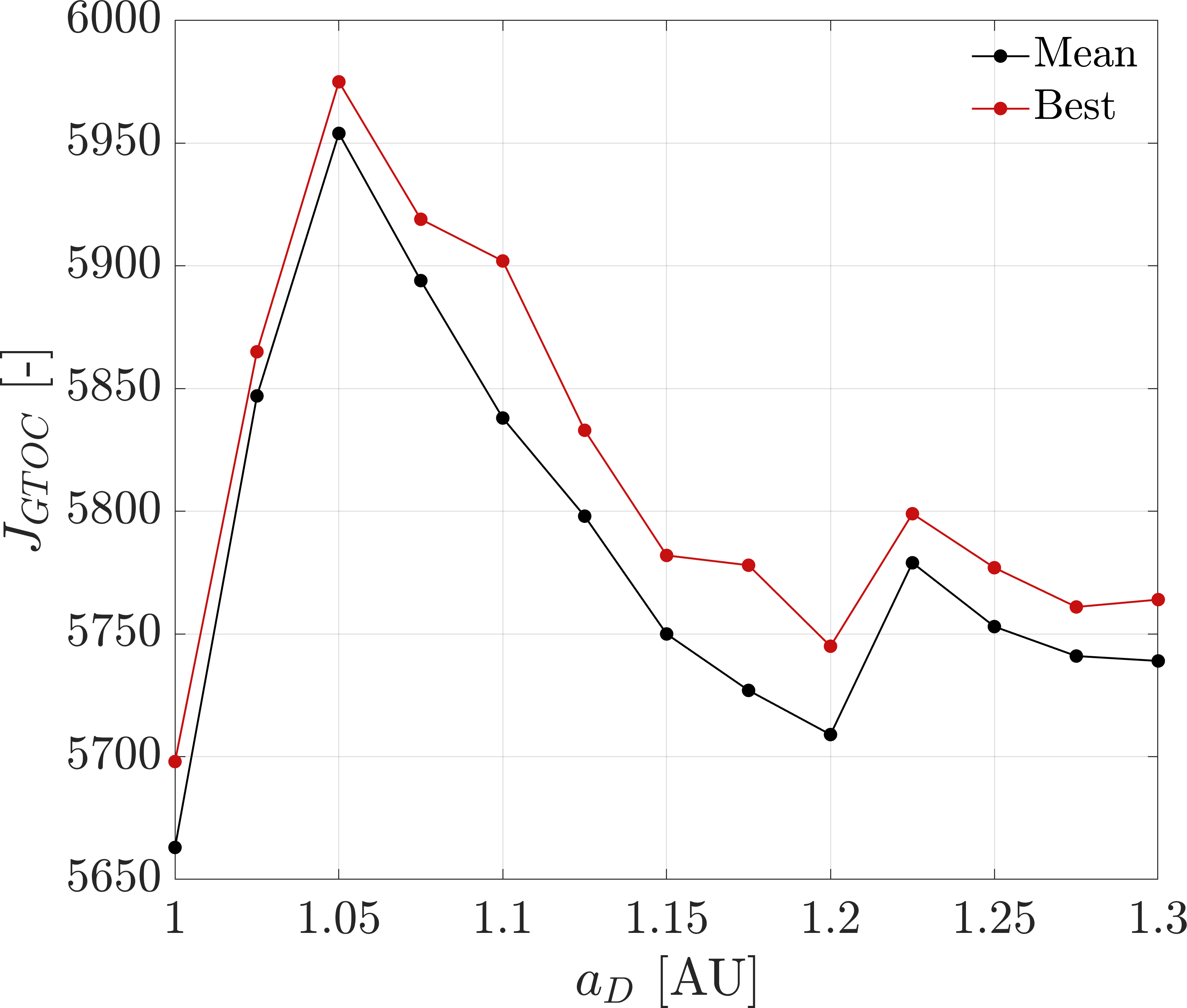}
\caption{Comparison of optimization results for different fixed ring radii. 5 runs for each radius were done, with the mean and best objectives shown.}
\label{fig:DispatcherRefinementConv1}
\end{center}
\end{figure}

{Different optimization algorithms were assessed and compared in the first stages of the competition, i.e. the MATLAB Surrogate Optimization Algorithm, Bayesian Optimization Algorithm, GA, PSO. 
Notably, the HPC was available for this task only about two days before the deadline. Until then, optimization runs were manually performed and repeated multiple times for promising solutions. In this case, a single run required 10 to 20 minutes (depending on settings) on a MacBook Pro, 2.3 GHz Quad-Core Intel Core i7. As a result, the final strategy to solve the optimization problem was to use the GA run on conventional laptops first to screen the solutions. Successively, a PSO algorithm was run for more extensive optimization of the promising solutions since it provided most stable solutions. For each case resulting from the analysis presented in Sec.~\ref{sec:TG}, five runs were performed on the Surrey's Space Centre's HPC cluster, providing the results of Fig.~\ref{fig:DispatcherRefinementConv1}. Each run required about 2hr of computations, for a population of 4,000 and stall iterations of 100.}

\section{The Refinement}\label{sec:5}
The refinement is the last step needed to obtain feasible solutions and to further reduce the motherships $\Delta v$.

\subsection{Asteroid Trajectory Refinement}
This section defines the asteroid trajectory refinement after they are allocated to the stations by the dispatcher. The time-optimal rendezvous problem is solved for each asteroid, taking the initial costates obtained from interpolation during the rendezvous table generation as initial guesses. The problem formulation and the dynamics of the time optimal rendezvous case are the same as those given in Sec.~\ref{toptfree}. {The Hamiltonian is the same as given in Eq. \ref{hamil}, however, for the rendezvous problem the shooting function now becomes }
\begin{equation}
   \boldsymbol{\phi}\left(\boldsymbol{\lambda}\left(t_{0}\right)\right)=\left[\boldsymbol{x}^{T}\left(t_{f}\right)-{\boldsymbol{x}_{f}}^{T}, H(t_f) - n_{D}{\lambda}_L(t_f),  \right]^{T}=\boldsymbol{0}, 
   \label{newphi}
\end{equation}
{in which $n_D$ is the mean motion of the Dyson ring.}
{Note that, by using the orbital representation of the states, the rendezvous constraint on the final Hamiltonian
\begin{equation}
H(t_f) -\boldsymbol{\lambda}(t_f)^T \dot{\boldsymbol{x}}_{f} = 0,    
\end{equation}
reduces to $H(t_f) - n_{D}{\lambda}_L(t_f) = 0$ as the true longitude is the only non-constant orbital element. Additionally, to account 
for differences in the final true anomaly that are multiples of $2\pi$, we have replaced the constraint $L(t_f) -L_f = 0$ with $\sin{(\frac{1}{2}(L(t_f) -L_f))} = 0$.}

\par
\begin{figure}
\begin{center}
\includegraphics[width=0.5\columnwidth]{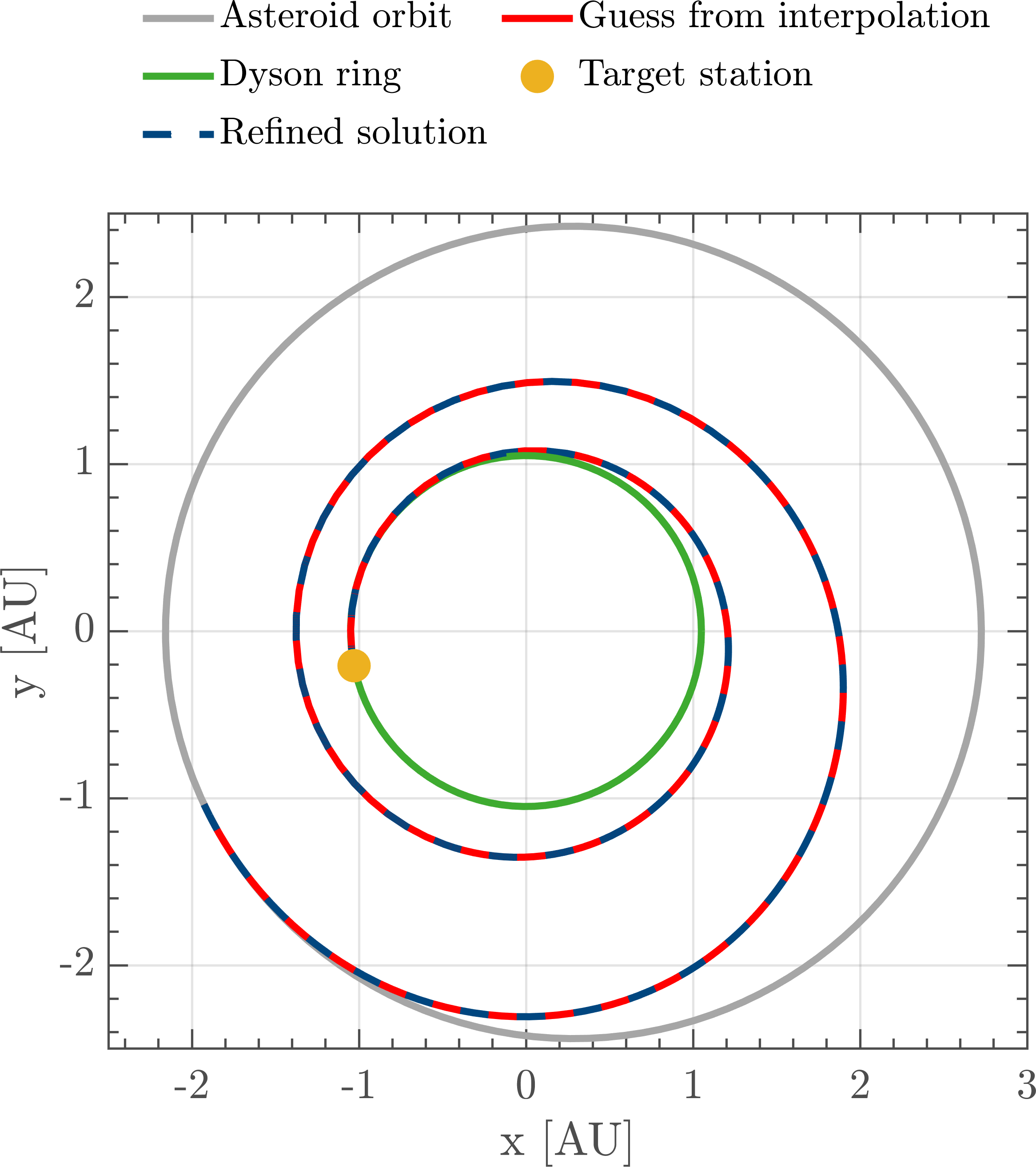}
\caption{Time optimal trajectory from initial guess vs refinement}
\label{fig:Rendezvous}
\end{center}
\end{figure}

 
Fig.~\ref{fig:Rendezvous} shows an example of an asteroid trajectory obtained by the interpolation and after the refinement. Clearly, the two trajectories are a close match.

The trajectory refinement was implemented in MATLAB and each trajectory was solved with simple shooting by the default algorithm implemented in \textit{fsolve}. Refining all transfers in the optimal solution required around 10 minutes on a MacBook Pro with a 2.3 GHz Quad-Core Intel Core i7 processor. 

\subsection{Deep Space Maneuvers}\label{ssec:DSM}
In the mothership chains described in Sec.~\ref{Sec:CAS} each transfer leg has only two $\Delta v$s, one at departure and one at arrival. However, the problem allows up to four $\Delta v$ maneuvers per leg. The last refinement introduced enables a single deep-space maneuver per leg, keeping all the encounter dates fixed. Initially, this approach was applied to re-optimize the legs containing asteroids not selected by the dispatcher {(these asteroids are always dropped and substituted with a DSM)}, but later, DSMs were enabled in all legs. 

An optimization problem was defined to determing the DSMs, considering the following optimization variables: $x$ is the fraction of the leg $\Delta t$, and $\boldsymbol{r}_{DSM}$ is the DSM position vector. The objective function to be minimized consists of the sum of the three $\Delta v$s computed by solving two Lamberts' arcs: the first one linking the departure asteroid to $\boldsymbol{r}_{DSM}$ in $x \Delta t$, and the second one linking  $\boldsymbol{r}_{DSM}$ to the arrival asteroid in $(1-x)\Delta t$. When an asteroid is removed from the chain, the guess for $x$ and $\boldsymbol{r}_{DSM}$ are related to the time and position corresponding to the skipped asteroid, which means that the DSM replaces the flyby. Instead, when no asteroids are removed, the guess for the maneuver time and position corresponds to half the leg.

For the optimal solution, the total $\Delta v$ per mothership after the introduction of the DSMs is reported in Table~\ref{table:MS_from_BS_test} and indicated as $\Delta v_{DSM}$. Although most maneuvers are small, the overall $\Delta v$ saving is about 1 km/s, most of it achieved in mothership 9, where an asteroid is dropped {(25 asteroids are visited instead of 26)}. Note that up to 22 DSMs are introduced in the mothership chains, almost one for each leg.

The optimization problem was solved with the interior point method implemented in MATLAB {\it fmincon}. The calculation of the DSMs for all the trajectories required around 15 minutes on a MacBook Pro with a 2.3 GHz Quad-Core Intel Core i7 processor.

\section{Optimal Solution}\label{sec:OS}

Table~\ref{table:MS_from_BS_test} illustrates the number of asteroids visited and the total $\Delta v$ spent by each of the 10 motherships. It can be seen that the number of asteroids visited remains in the 28-31 range, while the total $\Delta v$ of each mothership expenditure varies from 17.04 to 22.27 km/s. The total  $\Delta v$ spent is 194.90 km/s, which is taken into account in the performance index. \par
\begin{figure}[!ht]
\begin{center}
\includegraphics[width=.9\columnwidth]{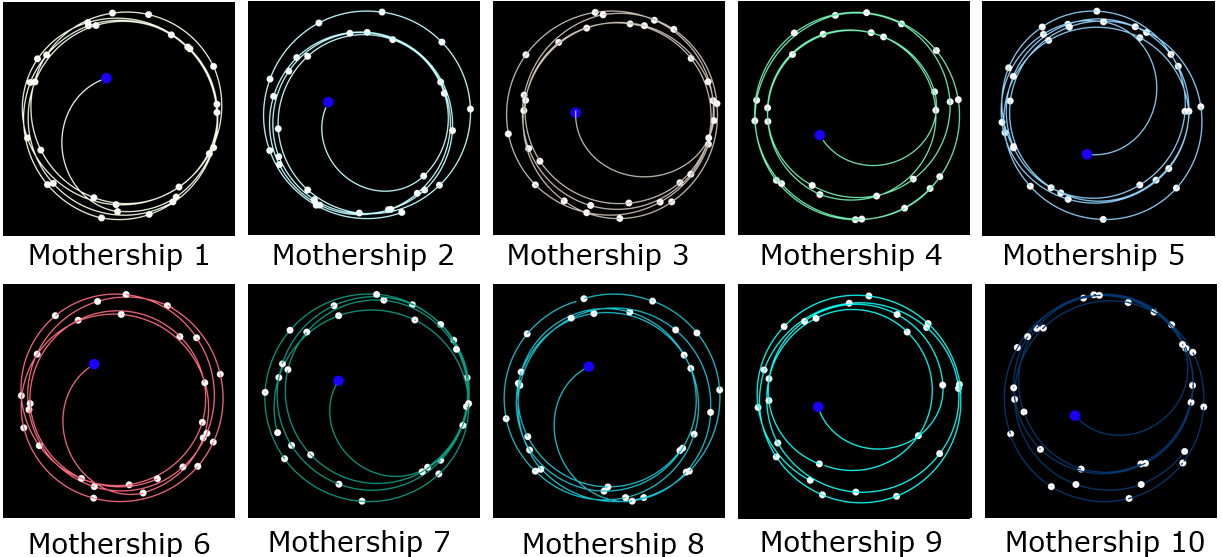}
\caption{Mothership Trajectories}
\label{fig:motherships}
\end{center}
\end{figure}

 
Fig.~\ref{fig:motherships} shows the trajectories taken by the motherships. They all start their journey from Earth and perform flyby maneuvers about the asteroids. 

\begin{table}[H]
\centering
\begin{tabular}{c|cccc|cc|cc}
\hline
Mothership &	\multicolumn{4}{c|}{Section~\ref{sec:2}}							&	\multicolumn{2}{c|}{Section~\ref{sec:4.1}}			&	\multicolumn{2}{c}{Section~\ref{sec:5}}			\\
	&	Ast	&	$\Delta v_{BS}$ & 		$\Delta v_{R1}$ &	$\Delta v_{R2}$ & Ast	&	$\Delta v_{trim}$	&	Ast	&	$\Delta v_{DSM}$ \\
	& Nr &	\multicolumn{3}{c|}{[km/s]} &	Nr	&	[km/s]		&	Nr	&	[km/s]		\\
\hline
1	&	30	&	22.80	&	19.89	&	19.43	&	29	&	19.16	&  29  &18.96			\\
2	&	32	&	21.09	&	17.67	&	17.46	&	31	&	17.22	& 31  &17.04			\\
3	&	32	&	25.19	&	21.65	&	20.68	&	32	&	20.68	&  32 &20.62			\\
4	&	30	&	22.35	&	19.64	&	19.39	&	28	&	18.85	&  28 &18.84		\\
5	&	31	&	25.65	&	21.73	&	21.3	&	31	&	21.30	&  31 &21.25		\\
6	&	32	&	26.3	&	22.59	&	22.25	&	28	&	18.78	&  28 &18.75		\\
7	&	31	&	21.92	&	18.94	&	18.59	&	29	&	18.49	&  29 &18.46		\\
8	&	32	&	27.84	&	24.03	&	23.26	&	30	&	22.28	&  30 &22.27		\\
9	&	31	&	24.18	&	21.08	&	20.26	&	26	&	18.47	&  25 &17.98		\\
10	&	32	&	25.98	&	22.46	&	21.75	&	30	&	20.69	&  30 &20.69		\\
\hline							
Total	&	313	&	243.32	&	209.67	&	204.39	&	294	&	195.91	&		293 & 194.90		\\
\hline
\end{tabular}
\caption{Mothership chains evolution in terms of asteroids and $\Delta v$, during the different phases of the optimization process.}
\label{table:MS_from_BS_test}
\end{table}

The final Dyson ring orbital elements are: $a_D = 1.05$ AU, $i_D= 0.85$ deg, $\Omega_D = 107.57$ deg. The initial phase of the first station was arbitrarily set to 0. Table~\ref{T:stations} illustrates the total number of asteroids going to each station, along with the total mass of each station once constructed. 
It can be seen that station 9 has the smallest mass upon construction ($\SI{1.2767e15}{\kilo\gram}$), and hence will be a factor in the performance index. 

\begin{table}[H]
\centering
\begin{tabular}{ccc}
\hline
Station no & No of Asteroids & Total mass ($\times 10^{15}$ kg) \\ \hline
1          & 27              & 1.2967 \\
2          & 26              & 1.2874 \\
3          & 29              & 1.2928 \\
4          & 28              & 1.2954 \\
5          & 24              & 1.3068\\
6          & 21              & 1.2860 \\
7          & 22              & 1.2768\\
8          & 24              & 1.2815\\
9          & 24              & 1.2767\\
10         & 23              & 1.2808 \\
11         & 21              & 1.2890 \\
12         & 24              & 1.2819 \\ \hline
\end{tabular}
\caption{Number of Asteroids and the total mass of each station upon construction}
\label{T:stations}
\end{table}

As our best solution was submitted at Nov 07, 2021 8:14 AM UTC, the performance index was obtained by evaluating Eq.~\eqref{eq:obj} 
\begin{equation}
    J_{GTOC}=B \cdot \frac{10^{-10} \cdot 1.2767\times 10^{15}}{1.05^{2} \sum_{k=1}^{10}\left(1 + \Delta v_i  / 50\right)^{2}} = 5992.3
\end{equation}
where the $\Delta v_i$ are taken from $\Delta v_{DSM}$ column of Table~\ref{table:MS_from_BS_test}. 

This was the solution with the third highest performance index by the deadline of GTOC 11. 

\section{Conclusions} \label{sec:C}
This paper presented the methodology developed by the team \lq\lq theAntipodes\rq\rq during the GTOC11 competition. Our approach allowed us to decouple the problem of motherships' chain generation and allocation of asteroids to stations. The interpolation approach adopted for the rendezvous table generation was extremely powerful as it enabled the estimation of all possible asteroids transfers to the Dyson ring at high accuracy. This allowed the dispatcher to work with a reliable value for the target function and made the final refinement an easy task.  

All the algorithms used to solve the problem {rely on light computations. Only at the very beginning of the competition and towards the end of it, an HPC was used: 1) to assist the transcription of motherships' chain generation problem (see Sec.\ref{Sec:CAS}), and 2) to} run multiple instances of the dispatcher to investigate the effect of the Dyson ring size and mitigate the variability of the optimal solution due to a stochastic solver (see Sec.\ref{sec:4.1}).


The simplified assumptions made in the generation of the mothership chains (a limited number of time of flights) and the preliminary estimate of the asteroids arrival masses with Edelbaum's equation were probably two major limitations of our approach. 

{The interpolation method adopted for the rendezvous tables can be potentially used in minimum-time rendezvous trajectory optimizations, leveraging the much better convergence properties of time-optimal orbital transfer with respect to the rendezvous ones. Additionally, solving the fixed acceleration problem made us realize that the Lagrangian multiplier of the mass can always be ignored in time-optimal transfers. This aspect will be discussed in a separate work currently in preparation. The method for balancing the mass distribution among the stations can be easily adopted for the design of multi active debris removal missions requiring a balance of removed mass among the missions. }


\section*{Acknowledgment}
The team is grateful to the University of Surrey for giving us access to their HPC.

\bibliographystyle{elsarticle-num} 
\bibliography{reference}

\begin{thebibliography}{10}
\expandafter\ifx\csname url\endcsname\relax
  \def\url#1{\texttt{#1}}\fi
\expandafter\ifx\csname urlprefix\endcsname\relax\def\urlprefix{URL }\fi
\expandafter\ifx\csname href\endcsname\relax
  \def\href#1#2{#2} \def\path#1{#1}\fi

\bibitem{mitchell1998introduction}
M.~Mitchell, An introduction to genetic algorithms, MIT press, 1998.

\bibitem{izzo2015revisiting}
D.~Izzo, Revisiting lambert’s problem, Celestial Mechanics and Dynamical
  Astronomy 121~(1) (2015) 1--15.

\bibitem{Kluever2011}
C.~A. Kluever, \href{http://arc.aiaa.org/doi/10.2514/1.51024}{{Using Edelbaum's
  Method to Compute Low-Thrust Transfers with Earth-Shadow Eclipses}}, Journal
  of Guidance, Control, and Dynamics 34~(1) (2011) 300--303.
\newblock \href {https://doi.org/10.2514/1.51024} {\path{doi:10.2514/1.51024}}.
\newline\urlprefix\url{http://arc.aiaa.org/doi/10.2514/1.51024}

\bibitem{rasotto2016multi}
M.~Rasotto, R.~Armellin, P.~Di~Lizia, Multi-step optimization strategy for
  fuel-optimal orbital transfer of low-thrust spacecraft, Engineering
  Optimization 48~(3) (2016) 519--542.

\bibitem{Eberhart1995}
R.~Eberhart, J.~Kennedy, {A new optimizer using particle swarm theory}, in:
  MHS'95. Proceedings of the Sixth International Symposium on Micro Machine and
  Human Science, Ieee, 1995, pp. 39--43.

\bibitem{schlueter2012nonlinear}
M.~Schlueter, Nonlinear mixed integer based optimization technique for space
  applications, Ph.D. thesis, University of Birmingham (2012).

\bibitem{Alemany2007}
K.~Alemany, R.~D. Braun, {Survey of global optimization methods for low-thrust,
  multiple asteroid tour missions} (2007).

\bibitem{Shapiro1992}
S.~C. Shapiro, {Encyclopedia of artificial intelligence second edition}, John,
  1992.

\bibitem{Wilt2010}
C.~M. Wilt, J.~T. Thayer, W.~Ruml, {A comparison of greedy search algorithms},
  in: Third Annual Symposium on Combinatorial Search, 2010.

\bibitem{Englander2012a}
J.~A. Englander, B.~A. Conway, T.~Williams, {Automated mission planning via
  evolutionary algorithms}, Journal of Guidance, Control, and Dynamics 35~(6)
  (2012) 1878--1887.

\bibitem{edelbaum}
T.~N. Edelbaum, Propulsion requirements for controllable satellites, ARS
  Journal 31 (1961) 1079--1089.
\newblock \href {https://doi.org/10.2514/8.5723} {\path{doi:10.2514/8.5723}}.

\bibitem{Bellome2021}
A.~Bellome, J.~P. {Sanchez Cuartielles}, J.~{Del Ser}, S.~Kemble, L.~Felicetti,
  {Efficiency of tree-search like heuristics to solve complex mixed-integer
  programming problems applied to the design of optimal space trajectories},
  Proceedings of the International Astronautical Congress, IAC (2021).

\bibitem{A&B}
J.~T. Betts, Very low-thrust trajectory optimization using a direct sqp method,
  Journal of Computational and Applied Mathematics 120~(1) (2000) 27--40.
\newblock \href {https://doi.org/https://doi.org/10.1016/S0377-0427(00)00301-0}
  {\path{doi:https://doi.org/10.1016/S0377-0427(00)00301-0}}.

\bibitem{diwu}
D.~Wu, F.~Jiang, J.~Li, Warm start for low-thrust trajectory optimization via
  switched system, Journal of Guidance, Control and Dynamics 47 (2021) 8--11.
\newblock \href {https://doi.org/10.2514/1.G005983}
  {\path{doi:10.2514/1.G005983}}.

\bibitem{gough2009gnu}
B.~Gough, GNU scientific library reference manual, Network Theory Ltd., 2009.

\bibitem{rasotto2016differential}
M.~Rasotto, A.~Morselli, A.~Wittig, M.~Massari, P.~Di~Lizia, R.~Armellin,
  C.~Valles, G.~Ortega, Differential algebra space toolbox for nonlinear
  uncertainty propagation in space dynamics, in: 6th International Conference
  on Astrodynamics Tools and Techniques (ICATT), 2016.

\end{thebibliography}

\end{document}